\documentclass{article}[12 pt]
\usepackage{csquotes}
\usepackage{amssymb}
\usepackage{amsmath}
\usepackage{amsthm}
\usepackage{graphicx}
\graphicspath{ {images/} }
\usepackage{mwe}
\usepackage{tikz} 
\usepackage{framed}
\usepackage{stackengine}
\usepackage[english]{babel}
\usepackage{empheq}
\usepackage{caption,subcaption}
\usepackage{amsmath,amssymb,tikz-cd}
\graphicspath{ {./images/} }
\usepackage[
backend=biber,
style=numeric,
]{biblatex}
\theoremstyle{definition}
\newtheorem{definition}{Definition}[section]

\theoremstyle{remark}

\title{The LG Fibration}
\author{\textbf{Daniel Livschitz}, \textbf{Weiqing Gu}\\ Claremont Graduate University,\\ Institute of Mathematical Sciences\\ Harvey Mudd College\\
\\$[$\textit{gu, dlivschitz}$]$@\textit{g.hmc.edu}}
\begin{document}
\maketitle
\begin{abstract}

 Deep Learning has significantly impacted the application of data-to-decision throughout research and industry, however, many Deep Learning methods lack a rigorous mathematical foundation, which creates situations where algorithmic results fail to be practically invertible. In this paper we present a nearly invertible mapping between $\mathbb{R}^{2^n}$ and $\mathbb{R}^{n+1}$ via a topological connection between $S^{2^n-1}$ and $S^n$. Throughout the paper we utilize the algebra of Multicomplex rotation groups and polyspherical coordinates to define two maps: the first is a contraction from $S^{2^n-1}$ to $\displaystyle \bigotimes^n_{k=1} SO(2)$, and the second is a projection from $\displaystyle \bigotimes^n_{k=1} SO(2)$ to $S^{n}$, which together form a composite map with an almost fibration structure that we call the LG Fibration. In analogy to the generation of Hopf Fibration using Hypercomplex geometry from $S^{(2n-1)} \mapsto CP^n$, our fibration uses Multicomplex geometry to project $S^{2^n-1}$ onto $S^n$. Additionally, we investigate the algebraic properties of the LG Fibration, ultimately deriving a distance difference function to determine which pairs of vectors have an invariant inner product under the transformation. The LG Fibration has applications to Machine Learning and Artificial Intelligence, in analogy to the current applications of Hopf Fibrations in adaptive UAV control. Furthermore, the ability to invert the LG Fibration for nearly all elements allows for the development of Machine Learning algorithms that may avoid the issues of uncertainty and reproducibility that currently plague contemporary methods.\cite{nagarajan2018impact} The primary result of this paper is a novel method of nearly invertible geometric dimensional reduction from $S^{2^n-1}$ to $S^n$, which has applications for further research in both mathematics and Artificial Intelligence, including but not limited to the fields of homotopy groups of spheres, algebraic topology, machine learning, and algebraic biology. 

\end{abstract}
\begin{flushleft}
\section{Introduction}
\end{flushleft}
\begin{flushleft}
\quad Since the discovery of Quaternions by Olinde Rodrigues\cite{gray1980olinde} and Sir William Hamilton\cite{hamilton1844ii} independently in 1840 and 1843, the anti-commutative Division Algebras of $\mathbb{H}$ and subsequently $\mathbb{O}$ have been the de facto generalization of the Complex Numbers in mathematical research. In 1877, Ferdinand Frobenius proved that all Real finite-dimensional associative division algebras must be isomorphic with either $\mathbb{R}$, $\mathbb{C}$, or $\mathbb{H}$. By 1923, Hurwitz's work on finite-dimensional unital real algebras was posthumously published, and proved that only in dimensions 1, 2, 4, and 8 are the conditions outlined in the Hurwitz Problem satisfied. These developments in representation theory, topology, and algebra set the stage for Heinz Hopf to publish his famous ``Hopf Fibration" in 1931. This established a many to one map from $S^3$ to $S^2$ using $S^1$ fibrations. Adam's Theorem further generalized the mapping to $S^7$ and $S^{15}$ using fibrations of $S^3$ and $S^7$ respectively, proving definitively that a Homotopy of spheres with a Hopf Invariant of 1 can only be found in dimensions 1, 2, 4, and 8. The algebras of $\mathbb{R}$, $\mathbb{C}$, $\mathbb{H}$, and $\mathbb{O}$ have a pronounced history in the formation of modern algebra, however their properties are only a part of the greater connection between generalizations of the Complex Numbers and Homotopy Groups of Spheres in high dimensions.
\end{flushleft}
\begin{flushleft}
\quad The Multicomplex Numbers are a $2^n$ dimensional commutative ring that extend the matrix representation of Complex Numbers to all higher orders of n, and contain a large number of geometric subgroups. Due to advancements in modern algebra and number theory, the once untenable zero divisors that plagued Multicomplex Numbers can now be explicitly defined with relative ease and separated from the set to form Abelian high dimensional multiplicative groups. These properties and more have allowed the Bicomplex numbers to find numerous applications throughout mathematical physics\cite{rochon2004bicomplex} and algebraic biology\cite{petoukhov2012dyadic}. The construction was first studied by Corrado Segre\cite{segre1892rappresentazioni}, who researched a commutative 4 dimensional non division ring extension of the Complex Numbers, now known as the Bicomplex Numbers. This established a method of generating consecutively higher dimensional commutative ring extensions of the Bicomplex Numbers. In 1991, Griffith Price formalized the general algebraic properties of Multicomplex Numbers in his seminal work ``An Introduction to Multicomplex Spaces and Functions"\cite{price2018introduction}. In this paper, we will use the abelian rotation groups in $\mathbb{C}_n$ to represent $\displaystyle \bigotimes^n_{k=1} SO(2)$, which serves as the link between $S^{2^n-1}$ and $S^n$. 
\end{flushleft}
\begin{flushleft}
\quad We begin section 2 with a review of the Hopf Fibration construction as well as additional background into their use in $\mathbb{R}^3$ spatial rotations. In sections 2.2- 2.3, we explore the polyspherical coordinate representation of a hypersphere as well as define the algebraic properties of the Bicomplex numbers. Throughout section 3 we explore the motivating example of the LG Fibration from $S^3$ to $S^2$, starting in 3.1 with the contraction from $S^3$ onto $\displaystyle \bigotimes^2_{k=1} SO(2)$. In 3.2 we define the projection map from $\displaystyle \bigotimes^2_{k=1} SO(2)$ onto $S^2$, and define the LG Fibration as the composition of the two mappings. In section 3.3 we explore how to visualize the LG Fibration through an intermediary mapping between the contraction and the projection. The LG Fibration is fully generalized in section 4 as we define the contraction from $S^{2^n-1}$ to $\displaystyle \bigotimes^n_{k=1} SO(2)$ as well as the second mapping that projects from $\displaystyle \bigotimes^n_{k=1} SO(2)$ to $S^n$. The definition of ``Almost Bijective Map" is given in section 4.4 alongside worked examples of the surjective vs bijective regions on $S^3$ and $S^7$. Section 5 explores the applications of the LG Fibration to Machine Learning by providing an explicit definition for the "Difference Function", which can be used to compute vector pairs that have an invariant inner product under the LG Fibration. We conclude the paper in Section 6 with a collection of selected graphs that display a portion of the unique underlying symmetries present in the LG Fibration.
\end{flushleft}
\begin{flushleft}
\quad The goal of this paper is to demonstrate a novel method for mapping $S^{2^n-1}$ to $S^n$ that has applications in fields as diverse as representation theory, algebraic topology, machine learning, and algebraic biology. The computationally convenient properties of a commutative high dimensional algebra provide the mathematical foundation for Geometric Dimensionality Reduction algorithms, generalized N-bit Algebras, as well as methods of graphing and data visualization that rely on Bicomplex rather than Quaternionic algebra.
\end{flushleft}
\begin{flushleft}
\section{Background}
\end{flushleft}

\subsection{\centering The Hopf Fibration}

\begin{flushleft}
The Hopf Fibration was a breakthrough in the field of algebraic topology due to the discovery of a connection between $S^3$ and $S^2$ through circular fibers. By associating a great circle in $S^3$ with a point in $S^2$, a geometric mapping with connections to unit quaternions and applications to spatial rotations was formed. The Hopf Fibration can be defined as the following:
\end{flushleft}
\begin{definition}[Hopf Fibration]
Let $a^2+b^2+c^2+d^2=1$,  $a,b,c,d\in \mathbb{R}
\\
H:S^3\mapsto S^2, \quad H(a,b,c,d)=\Big(a^2+b^2-c^2-d^2,2(ad-bc),2(ac-bd)\Big)
\\
\smallskip$
The mapping from $S^3$ to $S^2$ can be clearly seen in polyspherical coordinates$
\\
\smallskip H\Big(\cos(\theta_1)\cos(\theta_3),\sin(\theta_1)\cos(\theta_3),\cos(\theta_2)\sin(\theta_3),\sin(\theta_2)\sin(\theta_3)\Big)
\\
=\Big(\cos(2\theta_3),-\sin(2\theta_3)\sin(\theta_1-\theta_2),-\sin(2\theta_3)\cos(\theta_1+\theta_2)\Big) $
\end{definition}
\begin{flushleft}
Quaternions give a simple way to encode this axis-angle representation in four numbers, and can be used to apply (calculate) the corresponding rotation to a position vector $(\mathrm{x}, \mathrm{y}, \mathrm{z})$, representing a
point relative to the origin in $\mathbf{R}^{3}$ 
Therefore, a rotation of angle $\theta$ around the axis defined by the unit vector
$$
\vec{u}=\left(u_{x}, u_{y}, u_{z}\right)=u_{x} \mathbf{i}+u_{y} \mathbf{j}+u_{z} \mathbf{k}
$$
can be represented by a quaternion. This can be done using an extension of Euler's formula:
$$
\mathbf{q}=e^{\frac{\theta}{2}\left(u_{x} \mathbf{i}+u_{y} \mathbf{j}+u_{z} \mathbf{k}\right)}=\cos \frac{\theta}{2}+\left(u_{x} \mathbf{i}+u_{y} \mathbf{j}+u_{z} \mathbf{k}\right) \sin \frac{\theta}{2}
$$
It can be shown that the desired rotation can be applied to an ordinary vector $\mathbf{p}=\left(p_{x}, p_{y}, p_{z}\right)=p_{x} \mathbf{i}+p_{y} \mathbf{j}+p_{z} \mathbf{k}$ in 3-dimensional space, considered as a quaternion with a real coordinate equal to zero, by evaluating the conjugation of $\mathbf{p}$ by $\mathbf{q}$ :
$$
\mathbf{p}^{\prime}=\mathbf{q} \mathbf{p} \mathbf{q}^{-1}
$$
Notice that if $\mathbf{p}$ is defined as a real number, then $\mathbf{p}$ would remain invariant under the conjugation by $\mathbf{q} \mathbf{p} \mathbf{q}^{-1}$. Recall, using the Hamilton product, for each unit quaternion we can use it to define a rotation in $R^3$ via conjugation. 
$\mathbf{p}^{\prime}=\left(p_{x}^{\prime}, p_{y}^{\prime}, p_{z}^{\prime}\right)$ is the new position vector of the point after the rotation.

\

In a programmatic implementation, the conjugation is achieved by constructing a quaternion whose vector part is $\mathbf{p}$ and real part equals zero, and then performing the quaternion multiplication. The vector part of the resulting quaternion is the desired vector $\mathbf{p}^{\prime}$.
The key is that use a unit quaternion number to represent a rotation in $\mathbb{R}^3$, it splits the $\mathbb{R}^4$ into $\mathbb{R}$ and $\mathbb{R}^3$, which is the real and imaginary part of $\mathbb{H}$. This is still true for any other higher $\mathbb{R}^n$ which split as into $\mathbb{R}$ and $\mathbb{R}^{n-1}$ including Cayley numbers or high dimensional Clifford algebra. The polyspherical coordinates we are going to introduce below allow us to split $S^{n-1}$ into $S^{n-k-1}(\cos(\theta))$ and $S^{k-1}(\sin(\theta))$. 
\end{flushleft}
\begin{center}
\subsection{Polyspherical Coordinates}
\end{center}
\begin{flushleft}
In order to understand the techniques derived in this paper, as well as to heuristically follow the calculations performed, it is essential to discuss a particular method of generating higher dimensional spherical coordinates called polyspherical coordinates.
\end{flushleft}
\begin{flushleft}
The polyspherical coordinate construction states that the unit sphere in $\mathbb{R}^{n}$ can be generated inductively using polyspherical coordinates for $\mathbb{R}^{k}$ and $\mathbb{R}^{n-k}$. The relationship can be stated as the following:
\end{flushleft}
\begin{definition}[Polyspherical Coordinates]
Let $w\in S^{n-k-1},z\in S^{k-1},
\\\theta\in[0,2\pi)$
\end{definition}
\begin{center}
$S^{n-1}\equiv\big(w\cos(\theta),z\sin(\theta)\big)$
\end{center}
\begin{flushleft}
A special case of this construction is when we generalize from $\mathbb{R}^{n}$ to $\mathbb{R}^{n+1}$ by factoring in a $\cos(\theta_{n})$ to the previous polyspherical coordinates and adding an independent $\sin(\theta_{n})$\cite{vilenkin1993representation}:
\end{flushleft}
\begin{center}
$ S^{n-1}=\big\{\displaystyle\Big(\prod_{k=1}^{n-1} \cos \left(\theta_{k}\right)\Big)e_0+\sum_{k=1}^{n-1} \Big(\sin \left(\theta_{k}\right) \prod_{\ell=k+1}^{n-1} \cos \left(\theta_{\ell}\right) \Big)e_{k}\big| $
$\theta_k\in [0,\pi), 1\leq k< n-1,\theta_{n-1}\in [0,2\pi)\big\}$
\end{center}
\begin{center}
$ S^{n}\nolinebreak=\big\{\displaystyle\bigg(\Big(\prod_{k=1}^{n-1} \cos \left(\theta_{k}\right)\Big)e_0+\sum_{k=1}^{n-1} \Big(\sin \left(\theta_{k}\right) \prod_{\ell=k+1}^{n-1} \cos \left(\theta_{\ell}\right) \Big)e_{k}\bigg)\cos(\theta_{n})$  
\end{center}
\begin{center}
$+\sin(\theta_{n})e_{n}\big|\quad\theta_k\in [0,\pi), 1\leq k< n,\theta_{n}\in [0,2\pi)\big\}$
\end{center}
\begin{center}
$ =\big\{\displaystyle\Big(\prod_{k=1}^{n} \cos \left(\theta_{k}\right)\Big)e_0+\sum_{k=1}^{n} \Big(\sin \left(\theta_{k}\right) \prod_{\ell=k+1}^{n} \cos \left(\theta_{\ell}\right) \Big)e_{k}\big|
$
\end{center}
\begin{center}
$
\theta_k\in [0,\pi), 1\leq k< n,\theta_{n}\in [0,2\pi)\big\}$\nolinebreak
\end{center}
\begin{flushleft}
This form of $S^n$ will allow us to bridge the connection between Multicomplex Numbers and polyspherical coordinates. We will now focus on understanding the mathematics behind the product of rotations in the Multicomplex Algebra.
\end{flushleft}
\begin{center}
\subsection{Rotations In Bicomplex Space}
\end{center}
\begin{flushleft}
The Bicomplex numbers contain two unique complex terms i and j such that each squares to -1. We can use these elements to generate the full basis of the Bicomplex numbers by first taking the Power Set of the two basis, and then interpreting set inclusion as multiplication as well as the nullset as 1:
\end{flushleft}
\begin{center}
$\mathbb{G}_2\equiv\big\{i,j|$   $ i^2=j^2=-1, ij=ji\big\}$
\end{center}
\begin{center}
$\mathbb{P}(\mathbb{G}_2)=\Big\{\{\emptyset\},\{i\},\{j\},\{i,j\}\Big\}=\Big\{1,i,j,ij\Big\}$
\end{center}
\begin{flushleft}
We can now take the linear combination of the basis set over the Real numbers, along with the algebraic symmetries of the basis elements, and formally define a closed form expression for the Bicomplex numbers\cite{price2018introduction}:
\end{flushleft}
\begin{center}
$\mathbb{C}_2=\{a+bi+cj+dij|$   $a,b,c,d\in\mathbb{R}, 1,i,j,ij\in\mathbb{P}(\mathbb{G}_2)\}$
\end{center}
\begin{flushleft}
The Bicomplex numbers form a commutative ring, which means that the multiplication operation for elements in the set is well defined and Abelian. Furthermore, as a Banach Algebra, maps such as the exponential function are well defined in $\mathbb{C}_2$\cite{price2018introduction}. 
\end{flushleft}
\begin{center}
$(\mathbb{C}_2,+,\odot$ $,\|\|,\times)$
\end{center}
\begin{flushleft}
If we take the product of two independent rotations over the Bicomplex Numbers, the resulting space is a closed set that lies on $S^3$\cite{karakucs2015generalized}:\end{flushleft}

\begin{definition}[Simple $S^3$ Rotation Group]
Let $i,j\in \mathbb{G}_2, \theta\in[0,2\pi),\\ \psi\in[0,\pi)$
\begin{center}
$\mathbb{S}^3_{\mathbb{G}_2}\equiv\big\{\displaystyle e^{\displaystyle i\theta}\cdot e^{\displaystyle j\psi}\big\}$
$=\big\{\cos(\theta)\cos(\psi)+i\sin(\theta)\cos(\psi)+j\cos(\theta)\sin(\psi)+ij\sin(\theta)\sin(\psi)\big\}$
\end{center}
\end{definition}

\begin{flushleft}
We can verify that this surface is confined to $S^3$ by taking the inner product of any possible vector in $\mathbb{S}^3_{\mathbb{G}_2}$ with itself:
\end{flushleft}

\begin{center}
$w=\displaystyle e^{\displaystyle i\theta}\cdot e^{\displaystyle j\psi}$
\end{center}
\begin{center}
$\langle w,w\rangle=\cos^2(\theta)\cos^2(\psi)+\sin^2(\theta)\cos^2(\psi)+\cos^2(\theta)\sin^2(\psi)+\sin^2(\theta)\sin^2(\psi)
$
\end{center}
\begin{center}
$
=(\cos^2(\theta)+\sin^2(\theta))(\cos^2(\psi)+\sin^2(\psi))=1$
\end{center}
\begin{flushleft}
This proves that the magnitude squared is equal to 1, therefore:
\end{flushleft}
\begin{center}
$|w|=1$
\end{center}
\begin{flushleft}
Another method of factoring the product of two rotations gives insight into how we generate our projection:
\end{flushleft}
\begin{center}
$\displaystyle e^{\displaystyle i\theta}\cdot e^{\displaystyle j\psi}=(\cos(\theta)+i\sin(\theta))(\cos(\psi)+j\sin(\psi))$
\end{center}
\begin{center}
$=(\cos(\theta)+i\sin(\theta))\cos(\psi)+j(\cos(\theta)+i\sin(\theta))\sin(\psi)$
\end{center}
\begin{center}
$=\cos(\theta)\cos(\psi)+i\sin(\theta)\cos(\psi)+je^{\displaystyle i\theta}\sin(\psi)$
\end{center}
\section{LG Fibration in $\mathbb{R}^4$}

\subsection{\centering Mapping $S^3$ to $\mathbb{S}^3_{\mathbb{G}_2}$}

\begin{flushleft}
In this section we will explore the motivating example for the LG Fibration by constructing two maps that relate $S^{3}$ with $(S^1)^{\otimes 2}$, as well as $(S^1)^{\otimes 2}$ with $S^2$ through $\mathbb{T}^2$. The following diagram shows the relationships between the maps and spaces:
\end{flushleft}
\[
\begin{tikzcd}
S_p^{3} 
\arrow[r,twoheadrightarrow,"f"]&
\mathbb{S}^{3}_{\mathbb{G}_{2}} \arrow[r,leftrightarrow,"h"]
\arrow[dr,swap,twoheadrightarrow,"\mathbf{P}"]
&
\mathbb{T}^2_{2^{-1}}
\arrow[d,twoheadrightarrow,"\mu"]
\\
&&S^2
\end{tikzcd}
\]
\begin{center}
LG$(S_p^{3})=\mathbf{P}(f(S_p^{3}))=\mu(h(f(S_p^{3})))=S^2$
\end{center}
\begin{flushleft}
polyspherical coordinates allow us to map from $S^{2^n-1}$ onto the subspace $(S^1)^{\otimes n}$, while the commutative ring structure of Multicomplex Numbers allow us to represent $(S^1)^{\otimes n}$ as a rotation groups. We will define a map from $S^3$ onto $S^2$, as well as explore the LG Fibration mapping chain in order to illustrate the similarities and differences with the Hopf Fibration in low dimensions. The first step in defining the LG Fibration is establishing a particular coordinate system for $S^3$ to map onto  $\mathbb{S}^{3}_{\mathbb{G}_{2}}$:
\end{flushleft}
\begin{center}$
S_{p+}^{1}\equiv \{\big(\cos(\theta_1),\sin(\theta_1)\big)|$   $\theta_1\in [0,\pi)\}$
\end{center}
\begin{center}$
S_{p}^{1}\equiv \{\big(\cos(\theta_1),\sin(\theta_1)\big)|$   $\theta_1\in [0,2\pi)\}$
\end{center}
\begin{definition}[Particular $S^3$ Coordinate Representation] Let $\theta_1\in [0,\pi),\\ \theta_2\in [0,2\pi), \theta_3\in [0,\pi]$\\
\begin{center}$S_p^3\equiv \{\big(S_{p+}^{1}\cos(\theta_3),S_p^1\sin(\theta_3)\big)\}$
\end{center}
\begin{center}$
=\{\big(\cos(\theta_1)\cos(\theta_3),\sin(\theta_1)\cos(\theta_3),\cos(\theta_2)\sin(\theta_3),\sin(\theta_2)\sin(\theta_3)\big)\}$
\end{center}
\end{definition}
\begin{flushleft}
Now that we have a particular orientation for $S_p^{3}$, we can define a map that contracts $S_p^{3}$ to $\mathbb{S}^{3}_{\mathbb{G}_{2}}$ by grouping the input variables through their indices:
\end{flushleft}
\begin{definition}[Map From $S_p^{3}$ To $\mathbb{S}^{3}_{\mathbb{G}_{2}}$] Let $\theta_1\in [0,\pi),\theta_2\in [0,2\pi), \theta_3\in [0,\pi]$\\
\end{definition}
\begin{center}
$f(\theta_1,\theta_2,\theta_{3})=$
\end{center}
\begin{center}
$
\{\theta_1\mapsto \theta_{2},$ if $ \theta_{3}=\pi$ then $\theta_{2}\mapsto \theta_{2}\pm \pi, \theta_{3}\mapsto 0 \}$
\end{center}
\begin{flushleft}
Using Definition 3.1 and 3.2, we can evaluate explicitly how $S_p^{3}$ maps to $\mathbb{S}^{3}_{\mathbb{G}_{2}}$:
\end{flushleft}
\begin{center}
 $f\Big(\cos(\theta_1)\cos(\theta_3),\sin(\theta_1)\cos(\theta_3),\cos(\theta_2)\sin(\theta_3),\sin(\theta_2)\sin(\theta_3)\Big)$
\bigskip
$= \big(\cos(\theta_2)\cos(\theta_3),\sin(\theta_2)\cos(\theta_3),\cos(\theta_2)\sin(\theta_3),\sin(\theta_2)\sin(\theta_3)\big)$
\bigskip
$\cong \cos(\theta_2)\cos(\theta_3)+i\sin(\theta_2)\cos(\theta_3)+j\cos(\theta_2)\sin(\theta_3)+ij\sin(\theta_2)\sin(\theta_3)$
$=\displaystyle e^{i \theta_{2}}\cdot e^{j \theta_{3}}$
\end{center}
\begin{center}
$\theta_2\in [0,2\pi),\theta_3\in [0,\pi) $
\end{center}
\begin{flushleft}
This directly shows how the input of $S_p^{3}$ gets mapped onto $\mathbb{S}^{3}_{\mathbb{G}_{2}}$ through $f$:
\end{flushleft}
\begin{center}
$f(S_p^{3})=\mathbb{S}^{3}_{\mathbb{G}_{2}}$
\end{center}
\begin{center}
\subsection{The Projection From $\mathbb{S}^{3}_{\mathbb{G}_{2}}$ to $S^2$}
\end{center}
\begin{flushleft}
Unlike the Hopf Fibration, which provides a map for all components of the $S^3$ sphere, the LG Fibration preserves 2 of the components of $S^3$ in polyspherical coordinates and integrates the angles of rotation into the orientation of the reduced components. The factorization given at the end of section 2.4 provides us with the foundation for a projection from $\mathbb{S}^3_{\mathbb{G}_2}$ to $S^2$ by defining a map on the $je^{\displaystyle i\theta}\sin(\psi)$ component. Recall, we define $\mathbb{S}^3_{\mathbb{G}_2}$ as the following:
\end{flushleft}
\begin{center}
$\mathbb{S}^3_{\mathbb{G}_2}\equiv\big\{\displaystyle e^{\displaystyle i\theta}\cdot e^{\displaystyle j\psi}\big\}$
\\
$=\big\{\cos(\theta)\cos(\psi)+i\sin(\theta)\cos(\psi)+j\cos(\theta)\sin(\psi)+ij\sin(\theta)\sin(\psi)\big\}$
\\
$=\big\{\cos(\theta)\cos(\psi)+i\sin(\theta)\cos(\psi)+je^{\displaystyle i\theta}\sin(\psi)\big\}$
\end{center}
\begin{center}
$\theta_1\in [0,2\pi),\theta_2\in [0,\pi) $
\end{center}
\begin{center}
We can embed $\mathbb{S}^3_{\mathbb{G}_2}$ in $\mathbb{R}^4$ with the following coordinate representation:
\end{center}
\begin{center}
$x=\cos(\theta)\cos(\psi)$
\end{center}
\begin{center}
$y=\sin(\theta)\cos(\psi)$
\end{center}
\begin{center}
$z=\cos(\theta)\sin(\psi)$
\end{center}
\begin{center}
$w=\sin(\theta)\sin(\psi)$
\end{center}
\begin{center}
Similarly, we can see that polyspherical Coordinates for $S^2$ are defined as:
\end{center}
\begin{center}
$\bar{x}=\cos(\theta)\cos(\psi)$
\end{center}
\begin{center}
$\bar{y}=\sin(\theta)\cos(\psi)$
\end{center}
\begin{center}
$\bar{z}=\sin(\psi)$
\end{center}
\begin{center}
From this we can derive the following relationships:
\end{center}
\begin{center}
$\bar{x}=x$
\end{center}
\begin{center}
$\bar{y}=y$
\end{center}
\begin{center}
$\big|\bar{z}\big|=\sqrt{z^2+w^2}$
\end{center}
\begin{center}
We can now relate $\mathbb{S}^3_{\mathbb{G}_2}$ with $S^2$ through the following projection mapping:
\end{center}

\begin{center}
$\mathbf{P}\big(x+iy+jz+ijw\big)$
\end{center}
\begin{center}
$=\big(x,y,e^{i\big(\tan^{-1}(\frac{y}{x})-(\tan^{-1}(\frac{y}{x}) m o d \pi)\big)}\sqrt{z^2+w^2}\big)$
\end{center}
\begin{flushleft}
The motivation for this particular shift comes from the Bicomplex generalization of Euler's Identity that\cite{ozhan2022complex}:
\end{flushleft}
\begin{center}
$\displaystyle e^{i\pi}=e^{j\pi}=-1$
\end{center}
\begin{center}
If $\theta$ is within the range of $[0,\pi)$, then $\theta-(\theta \bmod \pi)=0$. 
\end{center}
\begin{center}
If $\theta$ is within the range of $[\pi,2\pi)$, then $\theta-(\theta \bmod \pi)=\pi$. 
\end{center}
\begin{flushleft}
We can use these modular relationships to construct the projection from $\mathbb{S}^{3}_{\mathbb{G}_{2}}$ to $ S^2$ by mapping $e^{i\theta}$ in $je^{i\theta}\sin(\psi)$ to $\pm 1$. 
\end{flushleft}
\begin{center}
$
e^{i \theta} \mapsto e^{i(\theta-(\theta \bmod \pi))}
$
\end{center}
\begin{flushleft}
This implies that the orientation of the z component is dependent on the angle of rotation $\theta$. We can normalize $\theta-(\theta \bmod \pi)$ by dividing by $\pi$ to indicate whether the result is 0 or 1. This produces in the following equality:
\end{flushleft}
\begin{center}
$\displaystyle e^{i\left(\theta-(\theta m o d \pi\right))}=(-1)^{\frac{\theta-(\theta m o d \pi)}{\pi}}$
\end{center}

\begin{flushleft}
This allows us to define the LG Fibration Projection as the map from $(S^1)^{\otimes2}$ embedded in $\mathbb{C}_2$ to $S^2$:
\end{flushleft}
\begin{definition}[Projection From $\mathbb{S}^{3}_{\mathbb{G}_{2}}$ To $ S^2 $] Let $\theta\in[0,2\pi),\psi\in[0,\pi)$

\end{definition}
\begin{center}
$\displaystyle \mathbf{P}(e^{\displaystyle i\theta}\cdot \smallskip e^{\displaystyle j\psi})=\mathbf{P}\big(\cos(\theta)\cos(\psi)+i\sin(\theta)\cos(\psi)+j\cos(\theta)\sin(\psi)+ij\sin(\theta)\sin(\psi)\big)
$
\end{center}
\begin{center}
$
=\mathbf{P}(\cos(\theta)\cos(\psi)+i\sin(\theta)\cos(\psi)+je^{\displaystyle i\theta}\sin(\psi))$
\end{center}
\begin{center}
$=\big(\cos(\theta)\cos(\psi),\sin(\theta)\cos(\psi),(-1)^{\frac{(\theta-(\theta m o d \pi))}{\pi}}\sin(\psi)\big)$
\end{center}
\begin{flushleft}
Using the LG Fibration Projection, we can recover a single cover of $S^2$ with half circle fibrations on the poles from $\mathbb{S}^3_{\mathbb{G}_2}$ as the following:
\end{flushleft}
\begin{center}
$x=\cos(\theta)\cos(\psi)$
\end{center}
\begin{center}
$y=\sin(\theta)\cos(\psi)$
\end{center}
\begin{center}
$z=(-1)^{\frac{(\theta-(\theta m o d \pi))}{\pi}}\sin(\psi)$
\end{center}
\begin{center}
$\theta\in [0,2\pi),\psi\in [0,\pi) $
\end{center}
\begin{center}
$\mathbf{P}(\mathbb{S}^3_{\mathbb{G}_2})=S^2$
\end{center}
\begin{center}
\subsection{Visualization of the LG Fibration in $\mathbb{R}^3$}
\end{center}
\begin{flushleft}
To visualize the points in $S^2$ that do not have a unique inverse to $\mathbb{S}^3_{\mathbb{G}_2}$, we must first project $\mathbb{S}^3_{\mathbb{G}_2}$ bijectively onto a surface in $\mathbb{R}^{3}$. $\mathbb{S}^3_{\mathbb{G}_2}$ is generated by the product of a full circle and of half a circle in $\mathbb{C}_2$, and is topologically one-to-one with $\frac{1}{2}$ of a torus in $\mathbb{R}^{3}$. This relationship is due to the construction of the torus as the product of two circles embedded in $\mathbb{R}^3$. In the given orientation, the half-circles on the northern and southern hemisphere of the torus are the fibers that map onto the northern and southern pole on $S^2$.
\end{flushleft}
\smallskip
\begin{center}
\begin{tabular}{|c|}
      \hline
      \includegraphics[width=.8\textwidth]{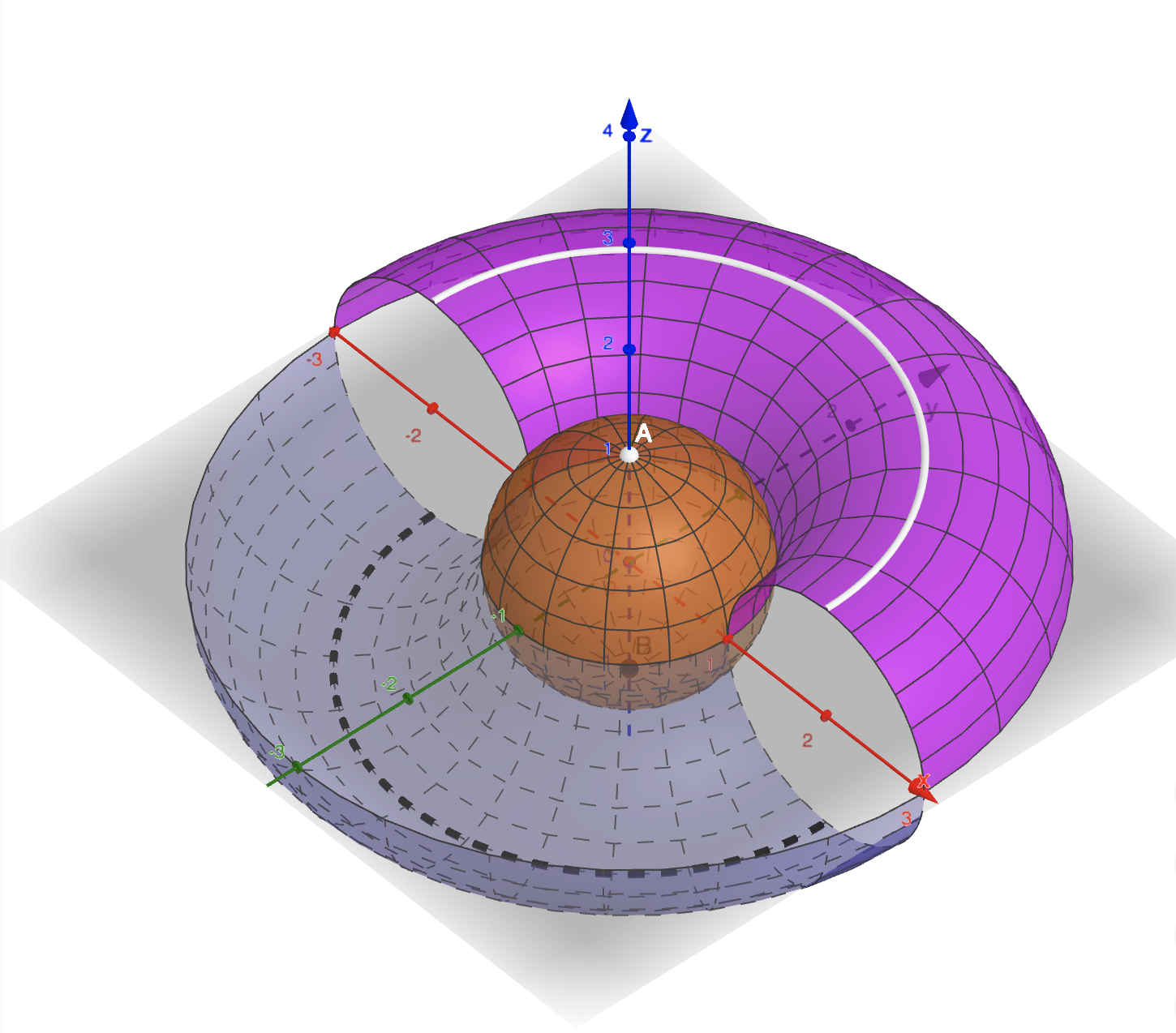}
      \\
\hline
\end{tabular}
\end{center}
\begin{center}
\end{center}
\begin{flushleft}
We can define the construction of the partial torus in $\mathbb{R}^3$ as the following:
\end{flushleft}
\begin{definition}[Half Torus in $\mathbb{R}^3$]
$\theta_1\in[0,2\pi), \theta_2\in[0,\pi), a_2\geq1 $
\end{definition}
\begin{center}
$\mathbb{T}^2_{2^{-1}}\equiv\big\{\Big(\cos(\theta_1)\big(a_2+\cos(\theta_2)\big),\sin(\theta_1)\big(a_2+\cos(\theta_2)\big),(-1)^{\frac{(\theta_1-(\theta_1 m o d \pi))}{\pi}}\sin(\theta_2)\Big)\big\}$
\end{center}
\begin{flushleft}
We can immediately show a one-to-one mapping between $\mathbb{S}^3_{\mathbb{G}_2}$ and $\mathbb{T}^2_{2^{-1}}$ by explicitly defining a map and as well as inverse mapping through their identical domain:
\end{flushleft}
\begin{definition}[Map from $\mathbb{S}^3_{\mathbb{G}_2}$ to Half Torus]
$\theta_1\in[0,2\pi), \theta_2\in[0,\pi) $
\end{definition}
\begin{center}
$h(e^{\displaystyle i\theta_1}\cdot e^{\displaystyle j\theta_2})=\Big(\cos(\theta_1)\big(a_2+\cos(\theta_2)\big),\sin(\theta_1)\big(a_2+\cos(\theta_2)\big),(-1)^{\frac{(\theta_1-(\theta_1 m o d \pi))}{\pi}}\sin(\theta_2)\Big)$
\end{center}
\begin{center}
$h(\mathbb{S}^3_{\mathbb{G}_2})=\mathbb{T}^2_{2^{-1}}$
\end{center}
\begin{flushleft}
All elements in $\mathbb{T}^2_{2^{-1}}$ and $\mathbb{S}^3_{\mathbb{G}_2}$ are uniquely defined by $\theta_1,\theta_2$, therefore we can simply define the inverse map as the following:
\end{flushleft}
\begin{center}
$h^{-1}\Big(\cos(\theta_1)\big(a_2+\cos(\theta_2)\big),\sin(\theta_1)\big(a_2+\cos(\theta_2)\big),(-1)^{\frac{(\theta_1-(\theta_1 m o d \pi))}{\pi}}\sin(\theta_2)\Big)
\newline
=e^{\displaystyle i\theta_1} \cdot e^{\displaystyle j\theta_2}$
\end{center}
\begin{center}
$\theta_1\in[0,2\pi), \theta_2\in[0,\pi) $
\end{center}
\begin{center}
$h^{-1}(\mathbb{T}^2_{2^{-1}})=\mathbb{S}^3_{\mathbb{G}_2}$
\end{center}
\begin{flushleft}
Now that we have defined a bijective map from $\mathbb{S}^3_{\mathbb{G}_2}$ to $\mathbb{T}^2_{2^{-1}}$, we must complete the relationship shown in the commutative diagram by defining a function $\mu(\mathbb{T}^2_{2^{-1}})$ that maps the partial torus onto $S^2$. This map can be seen as the limit of $\mathbb{T}^2_{2^{-1}}$ as $a_2\rightarrow 0.$
\end{flushleft}
\begin{definition}[Map from Half Torus to $S^2$]
$\theta_1\in[0,2\pi), \theta_2\in[0,\pi) $
\end{definition}
\begin{center}
$\mu\Big(\cos(\theta_1)\big(a_2+\cos(\theta_2)\big),\sin(\theta_1)\big(a_2+\cos(\theta_2)\big),(-1)^{\frac{(\theta_1-(\theta_1 m o d \pi))}{\pi}}\sin(\theta_2)\Big)$
\end{center}
\begin{center}$
=\big(\cos(\theta_1)\big(a_2+\cos(\theta_2)\big),\sin(\theta_1)\big(a_2+\cos(\theta_2)\big),(-1)^{\frac{(\theta_1-(\theta_1 m o d \pi))}{\pi}}\sin(\theta_2)\big)\newline
\lim a_2\rightarrow 0$
\end{center}
\begin{center}$
=\big(\cos(\theta_1)\cos(\theta_2),\sin(\theta_1)\cos\theta_2,(-1)^{\frac{(\theta_1-(\theta_1 m o d \pi))}{\pi}}\sin(\theta_2)\big)$
\end{center}
\begin{center}
$\mu(\mathbb{T}^2_{2^{-1}})=S^2$
\end{center}
\begin{flushleft}
Now that we have defined $S^3_p,\mathbb{S}^3_{\mathbb{G}_2},\mathbb{T}^2_{2^{-1}}$, and the maps between them, we are ready to define the LG Fibration for $S^3$. Lets recall the definition of the contraction and projection mappings:
\end{flushleft}
\begin{center}
$f\big(S^3_p\big)=\mathbb{S}^3_{\mathbb{G}_2}$
\end{center}
\begin{center}
$\mathbf{P}\big(\mathbb{S}^3_{\mathbb{G}_2}\big)= S^2$
\end{center}
\begin{flushleft}
These relationships allow us to decompose the LG Fibration of $S^3$ onto $S^2$ into two mappings, and can be formally defined as the following:
\end{flushleft}
\begin{definition}[LG Fibration for $S^3$]
\smallskip
Let $\theta_1\in[0,\pi),$ $\theta_2\in[0,2\pi),$ $ \theta_3\in[0,\pi]  $
\\
$LG:S_p^3\mapsto \mathbb{S}^3_{\mathbb{G}_2} \mapsto S^2$
\\
The contraction from $S_p^3$ to $\mathbb{S}^3_{\mathbb{G}_2}$ is achieved by mapping: $\theta_1\mapsto\theta_2$, and if $\theta_3=\pi$, then $\theta_3\mapsto 0$ and $\theta_2\mapsto\theta_2\pm \pi 
\\
$The projection from $\mathbb{S}^3_{\mathbb{G}_2}$ to $S^2$ is given by mapping: $je^{i\theta_2}\sin(\theta_3)\mapsto j(-1)^{\frac{\theta_2-(\theta_2 m o d \pi)}{\pi}}\sin(\theta_3)
\\
\\
\smallskip LG\Big(\cos(\theta_1)\cos(\theta_3),\sin(\theta_1)\cos(\theta_3),\cos(\theta_2)\sin(\theta_3),\sin(\theta_2)\sin(\theta_3)\Big)
\\
= \mathbf{P}\Big(\cos(\theta_2)\cos(\theta_3)+i\sin(\theta_2)\cos(\theta_3)+j\cos(\theta_2)\sin(\theta_3)+ij\sin(\theta_2)\sin(\theta_3)\Big)
\\
= \Big(\cos(\theta_2)\cos(\theta_3),\sin(\theta_2)\cos(\theta_3),(-1)^{\frac{\theta_2-(\theta_2 m o d \pi)}{\pi}}\sin(\theta_3)\Big)$
\end{definition}
\begin{center}
$LG\big(S^3_p\big)=\mathbf{P}\big(f\big(S^3_p\big)\big)= S^2$
\end{center}

\begin{flushleft}
\section{LG Fibration In High Dimensions}
\end{flushleft}

\subsection{\centering Rotations in Multicomplex Space}

\begin{flushleft}
Multicomplex numbers allow us to define an Abelian group embedded in $\mathbb{R}^{2^n}$, which is generated by the product of $n$ rotations. The benefit of a surface produced by $n$ orthogonal rotations in $\mathbb{C}_n$, is that the resulting set is confined to the surface of $S^{2^n-1}$ and only requires $n$ degrees of freedom to represent each point. By performing a modular reduction on the components of $(S^1)^{\otimes n}$, the LG Fibration allows us to reduce vectors from $\mathbb{R}^{2^n}$ to $\mathbb{R}^{n+1}$, as well as explicitly determine the orientation of the reduced components.
\end{flushleft}
\begin{flushleft}
We can generalize the construction of Multicomplex Numbers by defining the set of "simple" complex units that can be used with the Power set to generate the full basis for $\mathbb{C}_{n}:$  
\end{flushleft}
\begin{center}
$\mathbb{G}_{n}=\big\{i_k|$     $i^2_k=-1, i_mi_k=i_ki_m, 1\leq m,k\leq n\big\}$
\end{center}
\begin{center}
$\mathbb{P}\left(\mathbb{G}_{n}\right)\equiv\big\{\{\emptyset\},\{i_1\},...,\{i_1,i_2\},...,\{i_1,...,i_n\}\big\}
$
\end{center}
\begin{center}
$
=\big\{1,i_1,...,i_1i_2,...,i_1...i_n\big\}$
\end{center}
\begin{center}
$\mathbb{C}_{n}\equiv\big\{\displaystyle \sum_{u_k\in \mathbb{P}(\mathbb{G}_{n})}a_ku_k|$    $a_k\in \mathbb{R}, 1\leq k\leq 2^{n}\big\}$
\end{center}

\begin{flushleft}
Now that we have seen the basic tools required to multiply independent rotations together, we have the knowledge necessary to derive an explicit solution for each higher order LG Fibration.
\end{flushleft}
\begin{definition}[Simple $S^{2^n-1}$ Rotation Group]
Let $\theta_1\in [0,2\pi), \theta_m\in[0,\pi),
\\
1< j\leq n$
\end{definition}
\begin{center}
$\mathbb{S}^{2^{n}-1}_{\mathbb{G}_{n}}\equiv\big\{\displaystyle\prod^n_{k=1}e^{i_k\theta_k}|$   $i_k\in\mathbb{G}_{n}\}$
\\
$=\big\{\displaystyle\prod_{k=1}^{n} \cos \left(\theta_{k}\right)+\sum_{k=1}^{n}i_k\left(\sin \left(\theta_{k}\right) \prod_{\ell=k+1}^{n} \cos \left(\theta_{\ell}\right) \prod_{m=1}^{k-1} e^{i_m \theta_{m}}\right)|$   $i_k\in\mathbb{G}_{n}\big\}$
\end{center}
\begin{flushleft}
We inductively show the equality of the two definitions of $\mathbb{S}^{2^{n}-1}_{\mathbb{G}_{n}}$ in the following proof:
\end{flushleft}
\begin{center}
\subsection*{Proof of Equality}
\end{center}
\begin{center}
$\textbf{$2^1$ Dimensions}$
\end{center}
\begin{center}
If we set $n=1$, we can observe the following:
\end{center}
\begin{center}
$\displaystyle\prod^1_{k=1}e^{i_k\theta_k}=\cos(\theta_1)+i_1\sin(\theta_1)$
\end{center}
\begin{center}
Additionally, we can directly compute the second definition with $n=1$:
\end{center}
\begin{center}
$\displaystyle\prod_{k=1}^{1} \cos \left(\theta_{k}\right)+\sum_{k=1}^{1}i_k\left(\sin \left(\theta_{k}\right) \prod_{\ell=2}^{1} \cos \left(\theta_{\ell}\right) \prod_{m=1}^{0} e^{i_m \theta_{m}}\right)=\cos(\theta_1)+i_1\sin(\theta_1)$
\end{center}
\begin{center}
This allows us to conclude the proof of equality for $n=1$:
\end{center}
\begin{center}
$\displaystyle\prod^1_{k=1}e^{i_k\theta_k}=\prod_{k=1}^{1} \cos \left(\theta_{k}\right)+\sum_{k=1}^{1}i_k\left(\sin \left(\theta_{k}\right) \prod_{\ell=k+1}^{1} \cos \left(\theta_{\ell}\right) \prod_{m=1}^{k-1} e^{i_m \theta_{m}}\right)$
\end{center}
\
\begin{center}
$\textbf{$2^{n-1}$ Dimensions}$
\end{center}
\begin{center}
For our inductive step, we will assume that the following equality holds for $n-1$:
\end{center}
\begin{center}
$\displaystyle\prod^{n-1}_{k=1}e^{i_k\theta_k}\displaystyle=\prod_{k=1}^{n-1} \cos \left(\theta_{k}\right)+\sum_{k=1}^{n-1}i_k\left(\sin \left(\theta_{k}\right) \prod_{\ell=k+1}^{n-1} \cos \left(\theta_{\ell}\right) \prod_{m=1}^{k-1} e^{i_m \theta_{m}}\right)$
\end{center}
\begin{center}
$\textbf{$2^n$ Dimensions}$
\end{center}
\begin{center}
We must now show that the equality holds for $n$:
\end{center}
\begin{center}
$\displaystyle\prod^{n-1}_{k=1}e^{i_k\theta_k} \cdot\displaystyle e^{\displaystyle i_n \theta_{n}}=\displaystyle\prod^{n}_{k=1}e^{i_k\theta_k} $
\end{center}
\begin{center}
$\displaystyle\Big(\prod_{k=1}^{n-1} \cos\left(\theta_{k}\right)+\sum_{k=1}^{n-1}i_k\left(\sin \left(\theta_{k}\right) \prod_{\ell=k+1}^{n-1} \cos \left(\theta_{\ell}\right) \prod_{m=1}^{k-1} e^{i_m \theta_{m}}\right)\Big)\cdot \Big(\cos(\theta_n)+i_n\sin(\theta_n)\Big)$
\end{center}
\begin{center}
$\displaystyle=\Big(\prod_{k=1}^{n-1} \cos\left(\theta_{k}\right)+\sum_{k=1}^{n-1}i_k\left(\sin \left(\theta_{k}\right) \prod_{\ell=k+1}^{n-1} \cos \left(\theta_{\ell}\right) \prod_{m=1}^{k-1} e^{i_m \theta_{m}}\right)\Big)\cos \left(\theta_{n}\right)+\displaystyle \displaystyle\Big(\prod_{k=1}^{n-1} \cos\left(\theta_{k}\right)+\sum_{k=1}^{n-1}i_k\left(\sin \left(\theta_{k}\right) \prod_{\ell=k+1}^{n-1} \cos \left(\theta_{\ell}\right) \prod_{m=1}^{k-1} e^{i_m \theta_{m}}\right)\Big) i_n\sin(\theta_n)$
\end{center}
\begin{center}
$\displaystyle=\prod_{k=1}^{n} \cos\left(\theta_{k}\right)+\sum_{k=1}^{n-1}i_k\left(\sin \left(\theta_{k}\right) \prod_{\ell=k+1}^{n} \cos \left(\theta_{\ell}\right) \prod_{m=1}^{k-1} e^{i_m \theta_{m}}\right)+\displaystyle \Big(\prod_{m=1}^{n-1} e^{i_m \theta_{m}}\Big) i_n\sin(\theta_n)$
\end{center}
\begin{center}
$\displaystyle=\prod_{k=1}^{n} \cos\left(\theta_{k}\right)+\sum_{k=1}^{n}i_k\left(\sin \left(\theta_{k}\right) \prod_{\ell=k+1}^{n} \cos \left(\theta_{\ell}\right) \prod_{m=1}^{k-1} e^{i_m \theta_{m}}\right)$
\end{center}
\begin{center}
$ \boxed{ \displaystyle\prod^{n}_{k=1}e^{i_k\theta_k}=\prod_{k=1}^{n} \cos\left(\theta_{k}\right)+\sum_{k=1}^{n}i_k\left(\sin \left(\theta_{k}\right) \prod_{\ell=k+1}^{n} \cos \left(\theta_{\ell}\right) \prod_{m=1}^{k-1} e^{i_m \theta_{m}}\right)}$
\end{center}

\begin{center}
\subsection{Mapping $S^{2^n-1}$ to $\mathbb{S}^{2^{n}-1}_{\mathbb{G}_{n}}$}
\end{center}
\begin{flushleft}
The next step in generalizing the LG Fibration to higher dimensions is establishing a consistent orientation of $S^{2^n-1}$ to map onto $\mathbb{S}^{2^{n}-1}_{\mathbb{G}_{n}}$. We can begin by defining a method of generating a particular form of polyspherical Coordinates. In section 3.1, we defined the positive semi-circle $S_{p+}^{1}$ as the following:
\end{flushleft}
\begin{center}$
S_{p+}^{1}\equiv \{\big(\cos(\theta_1),\sin(\theta_1)\big)|$   $\theta_1\in [0,\pi)\}$
\end{center}
\begin{flushleft}
We then utilized the polyspherical structure to define a single cover of $S_p^3$ as the following:
\end{flushleft}
\begin{center}$
S_p^3\equiv \{\big(S_{p+}^{1}\cos(\theta_3),S_p^1\sin(\theta_3)\big)|$   $\theta_3\in [0,\pi]\}$
\end{center}
\begin{flushleft}
We can generalize the positive hemi-hysphere structure, which would allow us to define the single cover particular orientation of $S^{2^n-1}_{p+}$.
\end{flushleft}
\begin{center}
$S^{2^n-1}_{p+}\equiv\{\big(S^{2^{n-1}-1}_{p+}\cos(\theta_{2^n-1}),S^{2^{n-1}-1}_{p+}\sin(\theta_{2^n-1})\big)|$    $\theta_{2^n-1}\in [0,\pi)\}$
\end{center}
\begin{definition}[Particular Orientation of $S^{2^n-1}$]Let 
$\theta_{2^n-1}\in [0,\pi]$ 
\end{definition}
\begin{center}
$S^{2^n-1}_p\equiv\{\big(S^{2^{n-1}-1}_{p+}\cos(\theta_{2^n-1}),S^{2^{n-1}-1}_p\sin(\theta_{2^n-1})\big)\}$
\end{center}

\begin{center}
\subsection*{Example Map From $S^7$ to $\mathbb{S}^{7}_{\mathbb{G}_{3}}$}
\end{center}
\begin{flushleft}
Using Definition 4.1 and 4.2, we can define $S^{7}$ and $\mathbb{S}^{7}_{\mathbb{G}_{3}}$ as well as explictly show how $f$ maps from $S^{7}$ to $\mathbb{S}^{7}_{\mathbb{G}_{3}}$. This example allows us to demonstrate how the procedure generalizes for all $S^{2^n-1}\mapsto\mathbb{S}^{2^n-1}_{\mathbb{G}_{n}}$:
\end{flushleft}
\begin{center}$f(
\begin{bmatrix}
\cos(\theta_1)\cos(\theta_3)\cos(\theta_7)\\
\sin(\theta_1)\cos(\theta_3)\cos(\theta_7)\\
\cos(\theta_2)\sin(\theta_3)\cos(\theta_7)\\
\sin(\theta_2)\sin(\theta_3)\cos(\theta_7)\\
\cos(\theta_4)\cos(\theta_6)\sin(\theta_7)\\
\sin(\theta_4)\cos(\theta_6)\sin(\theta_7)\\
\cos(\theta_5)\sin(\theta_6)\sin(\theta_7)\\
\sin(\theta_5)\sin(\theta_6)\sin(\theta_7)
\end{bmatrix})=\begin{bmatrix}
\cos(\theta_5)\cos(\theta_3)\cos(\theta_7)\\
\sin(\theta_5)\cos(\theta_3)\cos(\theta_7)\\
\cos(\theta_5)\sin(\theta_3)\cos(\theta_7)\\
\sin(\theta_5)\sin(\theta_3)\cos(\theta_7)\\
\cos(\theta_5)\cos(\theta_3)\sin(\theta_7)\\
\sin(\theta_5)\cos(\theta_3)\sin(\theta_7)\\
\cos(\theta_5)\sin(\theta_3)\sin(\theta_7)\\
\sin(\theta_5)\sin(\theta_3)\sin(\theta_7)
\end{bmatrix}$
\end{center}
\begin{center}
$\cong \cos(\theta_5)\cos(\theta_3)\cos(\theta_7)+i_1\sin(\theta_5)\cos(\theta_3)\cos(\theta_7)+i_2\cos(\theta_5)\sin(\theta_3)\cos(\theta_7)+i_1i_2\sin(\theta_5)\sin(\theta_3)\cos(\theta_7)+i_3\cos(\theta_5)\cos(\theta_3)\sin(\theta_7)+i_1i_3\sin(\theta_5)\cos(\theta_3)\sin(\theta_7)+i_2i_3\cos(\theta_5)\sin(\theta_3)\sin(\theta_7)+i_1i_2i_3\sin(\theta_5)\sin(\theta_3)\sin(\theta_7)$
\end{center}
\begin{center}
$=\displaystyle e^{i_1 \theta_{5}}\cdot e^{i_2 \theta_{3}}\cdot e^{i_3 \theta_{7}}$
\end{center}
\begin{center}
$\theta_5\in [0,2\pi),\theta_3\in [0,\pi),\theta_7\in [0,\pi)$
\end{center}

\begin{flushleft}
Now that we have defined the particular orientation of $S^{2^n-1}$ that is in a similar form to $\mathbb{S}^{2^{n}-1}_{\mathbb{G}_{n}}$, we must explicitly define how to complete the map from $S_p^{2^n-1}\mapsto\mathbb{S}^{2^{n}-1}_{\mathbb{G}_{n}}$. This requires defining which set of indices for the input variables should be mapped to the same index in the output. Under our defined orientation of $S_p^{2^n-1}$, the input variables of $\theta_k$ for $S_{p+}^{2^{n-1}-1}$ are from $1\leq k\leq 2^{n-1}-1$, while the $\theta_m$ for $S_p^{2^{n-1}-1}$ are from $2^{n-1}\leq j\leq 2^{n}-2$. The surface $\mathbb{S}^{2^{n}-1}_{\mathbb{G}_{n}}$ only has $n$ degrees of freedom, so we must partition the $2^n-1$ input variables into n distinct sets:
\end{flushleft}
\begin{center}
$\Theta^{2^m-1}_0=\{2^m-1\}$
\end{center}
\begin{center}
$\Theta^{2^m-1}_{\ell}=\{k+2^{m+\ell-1}-1|$   $k\in\Theta^{2^m-1}_{\ell-1}\}\cup\Theta^{2^m-1}_{\ell-1}$
\end{center}
\begin{flushleft}
Finally, we can define the specific mappings from $S_p^{2^n-1}$ to $\mathbb{S}^{2^{n}-1}_{\mathbb{G}_{n}}$ using relations on the indices of the input variables:
\end{flushleft}
\begin{definition}[Map From $S_p^{2^n-1}$ to $\mathbb{S}^{2^{n}-1}_{\mathbb{G}_{n}}$]
\end{definition}
\begin{center}
$f(\theta_1,...,\theta_{2^n-1})=
 \Bigg\{ \begin{array}{rcl} \theta_k\mapsto \theta_{2^m-1} \quad\quad\quad\quad k\in \Theta^{2^m-1}_{n-m}, 1<m\leq n \\ 
 \theta_{\ell}\mapsto \theta_{2^n-n} \quad\quad\quad\quad \quad\quad\quad\quad\quad\quad \ell\in \Theta^{1}_{n-1}
  \\ 
\theta_{2^n-n}\mapsto \theta_{2^n-n}\pm \pi, \theta_{2^n-1}\mapsto 0  \quad \theta_{2^n-1}=\pi\end{array}$
\end{center}
\begin{center}
$f(S^{2^n-1})=\mathbb{S}^{2^n-1}_{\mathbb{G}_{n}}$
\end{center}

\begin{flushleft}
 The map $f$ can be heuristically interpreted as grouping vertical columns of input variables when $S^{2^n-1}_p$ is written in vector form. Using Definition 4.3 and the worked example of $S^7$, we can explicitly define $f:S^{2^n-1}_p\mapsto\mathbb{S}^{2^n-1}_{\mathbb{G}_{n}}$ as the following:
\end{flushleft}
\begin{center}$
f(S^{2^n-1}_p)=\displaystyle e^{i_1 \theta_{2^n-n}}\prod_{\ell=2}^{n} e^{i_{\ell} \theta_{2^{\ell}-1}}
$\end{center}
\begin{center}
$\theta_{2^n-n}\in [0,2\pi),$  $\theta_{2^{\ell}-1}\in [0,\pi), 1<\ell \leq n$
\end{center}
\begin{center}
\subsection{General LG Fibration Projection}
\end{center}
\begin{flushleft}
In section 3.1 we derived the following partial factorization for the product of $n$ rotations in $\mathbb{C}_n$:
\end{flushleft}
\begin{center}
$\displaystyle \displaystyle\prod_{k=1}^{n}e^{ i_k \theta_{k}}=\prod_{k=1}^{n} \cos\left(\theta_{k}\right)+\sum_{k=1}^{n}i_k\left(\sin \left(\theta_{k}\right) \prod_{\ell=k+1}^{n} \cos \left(\theta_{\ell}\right) \prod_{m=1}^{k-1} e^{i_m \theta_{m}}\right)$
\end{center}
\begin{flushleft}
We can now use this newly derived product form to generalize the LG Fibration Projection to map from $\mathbb{S}^{2^{n}-1}_{\mathbb{G}_{n}}$ to $S^n$ by making the following substitution:
\end{flushleft}
\begin{center}$
\displaystyle\prod_{m=1}^{k-1} e^{i_m \theta_{m}} \mapsto \prod_{m=1}^{k-1} e^{i_m(\theta_{m}-(\theta_{m} \bmod \pi))}
$\end{center}

\begin{flushleft}
By applying the generalized LG Fibration Projection to $\mathbb{S}^{2^{n}-1}_{\mathbb{G}_{n}}$, we recover the vector form of polyspherical coordinates in $\mathbb{R}^{n+1}$ with reduced components scaled by 1 or -1 depending on if $\frac{\theta_{k}-\left(\theta_{k} \bmod \pi\right)}{\pi}$ is 0 or 1:
\end{flushleft}
\begin{center}
$\mathbf{P}(\displaystyle\prod_{k=1}^{n} e^{\displaystyle i_k \theta_{k}})=\mathbf{P}(\displaystyle\prod_{k=1}^{n} \cos \left(\theta_{k}\right)+\sum_{k=1}^{n} i_{k} \sin \left(\theta_{k}\right) \prod_{\ell=k+1}^{n} \cos \left(\theta_{\ell}\right) \prod_{m=1}^{k-1} e^{i_{m} \theta_{m}}) $
\end{center}
\begin{flushleft}
We can now perform the LG Fibration Projection and map $1= e_0$ as well as $i_k= e_k$ for the vector representation in $\mathbb{R}^{n+1}$:
\end{flushleft}
\begin{center}
$= \displaystyle\bigg(\prod_{k=1}^{n} \cos \left(\theta_{k}\right)\bigg)e_0+\sum_{k=1}^{n} \bigg(\sin \left(\theta_{k}\right) \prod_{\ell=k+1}^{n} \cos \left(\theta_{\ell}\right) \prod_{m=1}^{k-1} e^{i_{m}\left(\theta_{m}-(\theta_{m} m o d \pi\right))}\bigg)e_{k} $
\end{center}
\begin{flushleft}
We can rewrite the final product of $e^{i_{m}\left(\theta_{m}-(\theta_{m} m o d \pi\right))}$ as the following:
\end{flushleft}
\begin{center}
$\displaystyle e^{i_{m}\left(\theta_{m}-(\theta_{m} m o d \pi\right))}=(-1)^{\frac{\theta_{m}-(\theta_{m} m o d \pi)}{\pi}}$
\end{center}
\begin{flushleft}
This brings us to the final form for the LG Fibration Projection:
\end{flushleft}
\begin{flushleft}
 \boxed{ \begin{array}{rcl} \mathbf{P}\big(\displaystyle\prod_{k=1}^{n}e^{\displaystyle i_k \theta_{k}}\big) = \displaystyle\bigg(\prod_{k=1}^{n} \cos \left(\theta_{k}\right)\bigg)e_0+\sum_{k=1}^{n} \bigg(\sin \left(\theta_{k}\right) \prod_{\ell=k+1}^{n} \cos \left(\theta_{\ell}\right) \prod_{m=1}^{k-1} (-1)^{\frac{(\theta_{m}-(\theta_{m} m o d \pi))}{\pi}}\bigg)e_{k}\end{array}}
 \end{flushleft}
 \begin{flushleft}
The result is the vector representation of $S^{2^n-1}$ in polyspherical Coordinates, thus we have proven that the LG Fibration Projection successfully maps from $\mathbb{S}^{2^{n}-1}_{\mathbb{G}_{n}}$ onto $S^n$, and can be seen as the following:
\end{flushleft}
\begin{center}
$\mathbf{P}\big(\mathbb{S}^{2^{n}-1}_{\mathbb{G}_{n}}\big)= S^{n}$
\end{center}
\begin{center}
\subsection{Almost Bijectivity}
\end{center}
\begin{flushleft}
We define an Almost Bijective Mapping as a map with a measure zero kernal. We show that the LG Fibration map is Almost Bijective by explictly defining the subsets for which the uniqueness of the inverse image is not preserved by the mapping. The generalized LG Fibration is defined as the following:
\end{flushleft}
\begin{center}
$\mathbf{P}(\displaystyle\prod_{k=1}^{n}e^{\displaystyle i_k \theta_{k}})\displaystyle=\bigg(\displaystyle\prod_{k=1}^{n} \cos \left(\theta_{k}\right)\bigg)e_0+\sum_{k=1}^{n}\bigg( \sin \left(\theta_{k}\right) \prod_{\ell=k+1}^{n} \cos \left(\theta_{\ell}\right) \prod_{m=1}^{k-1} (-1)^{\frac{(\theta_{m}-(\theta_{m} m o d \pi))}{\pi}}\bigg)e_k$
\end{center}
\begin{center}
$\theta_1\in [0,2\pi), \theta_k\in [0,\pi), 1<k\leq n$
\end{center}
\begin{flushleft}
If we restrict $\theta_k$ to $ [0,\pi)$, then necessarily $\displaystyle e^{i_{k}\left(\theta_{k}-(\theta_{k} m o d \pi\right))}=1$ for $1<k\leq n$. This allows us to simplify the mapping as the following:
\end{flushleft}
\begin{center}
$\mathbf{P}(\displaystyle\prod_{k=1}^{n}e^{\displaystyle i_k \theta_{k}})\displaystyle=\bigg(\displaystyle\prod_{k=1}^{n} \cos \left(\theta_{k}\right)\bigg)e_0+\bigg(\displaystyle\sin \left(\theta_{1}\right)\prod_{k=2}^{n} \cos \left(\theta_{k}\right)\bigg)e_1+\displaystyle (-1)^{\frac{(\theta_1-(\theta_1 m o d \pi))}{\pi}} \sum_{k=2}^{n}\bigg(\sin \left(\theta_{k}\right) \prod_{\ell=k+1}^{n} \cos \left(\theta_{\ell}\right) \bigg)e_k$
\end{center}
\begin{flushleft}
We can begin by working through the example of $\mathbb{S}^{3}_{\mathbb{G}_{2}}$ mapped Almost Bijectievly onto $S^2$, and understanding for which values the mapping is not unique:
\end{flushleft}
\begin{center}
$\displaystyle \mathbf{P}(e^{\displaystyle i_1\theta_1}\cdot  e^{\displaystyle i_2\theta_2})=
$
\end{center}
\begin{center}
$
\cos(\theta_1)\cos(\theta_2)e_0+\sin(\theta_1)\cos(\theta_2)e_1+(-1)^{\frac{(\theta_1-(\theta_1 m o d \pi))}{\pi}}\sin(\theta_2)e_2$
\end{center}
\begin{flushleft}
In this coordinate representation of the mapping, $\displaystyle e^{\displaystyle i_1\theta_1}$ acts as the horizontal counterclockwise rotation from the x to the y axis, and $e^{\displaystyle i_2\theta_2}$ acts as the vertical counterclockwise rotation from the xy plane to the z axis. We can observe the following for $\displaystyle \theta_2=\frac{\pi}{2}$:
\end{flushleft}
\begin{center}
$\displaystyle \mathbf{P}( e^{\displaystyle i_1\theta_1}\cdot \displaystyle e^{\displaystyle i_2 \frac{ \pi}{2}})=
$
\end{center}
\begin{center}
$
\cos(\theta_1)\cos(\frac{\pi}{2})e_0+\sin(\theta_1)\cos(\frac{\pi}{2})e_1+(-1)^{\frac{(\theta_1-(\theta_1 m o d \pi))}{\pi}}\sin(\frac{\pi}{2})e_2
$
\end{center}
\begin{center}
$=(-1)^{\frac{(\theta_1-(\theta_1 m o d \pi))}{\pi}}e_2$
\end{center}
\begin{flushleft}
We see that the following elements defined by $\mathbf{P}(e^{\displaystyle i_1\theta_1}\cdot \displaystyle e^{\displaystyle i_2 \frac{ \pi}{2}})$ map to the same point in $\mathbb{R}^3$:
\end{flushleft}
\begin{center}
$\mathbf{P}(e^{\displaystyle i_1\theta_1}\cdot \displaystyle e^{\displaystyle i_2 \frac{ \pi}{2}})= \left\{ \begin{array}{rcl} e_2 \quad\quad \theta_1\in [0,\pi) \\ - e_2  \quad \theta_1\in [\pi,2\pi)\end{array}\right.$   
\end{center}
\begin{flushleft}
We can extend the example to $\mathbb{S}^{7}_{\mathbb{G}_{3}}$ mapped Almost Bijectievly onto $S^3$:
\end{flushleft}

\begin{center}
$\displaystyle \mathbf{P}(e^{\displaystyle i_1\theta_1}\cdot e^{\displaystyle i_2\theta_2}\cdot e^{\displaystyle i_3\theta_3})=
$
\end{center}
\begin{center}
$
\cos(\theta_1)\cos(\theta_2)\cos(\theta_3)e_0
$
\end{center}
\begin{center}
$
+\sin(\theta_1)\cos(\theta_2)\cos(\theta_3)e_1
$
\end{center}
\begin{center}
$
+(-1)^{\frac{(\theta_1-(\theta_1 m o d \pi))}{\pi}}\sin(\theta_2)\cos(\theta_3)e_2
$
\end{center}
\begin{center}
$
+(-1)^{\frac{(\theta_1-(\theta_1 m o d \pi))}{\pi}}\sin(\theta_3)e_3$
\end{center}
\begin{flushleft}
We can observe a similar effect to the previous example for $\displaystyle \theta_3=\frac{\pi}{2}$:
\end{flushleft}
\begin{center}
$\displaystyle \mathbf{P}(e^{\displaystyle i_1\theta_1}\cdot e^{\displaystyle i_2\theta_2}\cdot e^{\displaystyle i_3\frac{\pi}{2}})=
$
\end{center}
\begin{center}
$
\cos(\theta_1)\cos(\theta_2)\cos(\frac{\pi}{2})e_0
$
\end{center}
\begin{center}
$
+\sin(\theta_1)\cos(\theta_2)\cos(\frac{\pi}{2})e_1
$
\end{center}
\begin{center}
$
+(-1)^{\frac{(\theta_1-(\theta_1 m o d \pi))}{\pi}}\sin(\theta_2)\cos(\frac{\pi}{2})e_2
$
\end{center}
\begin{center}
$
+(-1)^{\frac{(\theta_1-(\theta_1 m o d \pi))}{\pi}}\sin(\frac{\pi}{2})e_3$
\end{center}
\begin{center}
$=(-1)^{\frac{(\theta_1-(\theta_1 m o d \pi))}{\pi}}e_3$
\end{center}
\begin{center}
$\mathbf{P}(e^{\displaystyle i_1\theta_1}\cdot e^{\displaystyle i_2\theta_2}\cdot e^{\displaystyle i_3\frac{\pi}{2}})= \left\{ \begin{array}{rcl} e_3 \quad\quad \theta_1\in [0,\pi) \\ - e_3  \quad \theta_1\in [\pi,2\pi)\end{array}\right.$   
\end{center}
\begin{flushleft}
We can observe that if $\theta_k=\displaystyle\frac{\pi}{2}$, the components with a basis of $e_0...e_{k-1}$ are all scaled by $\cos(\frac{\pi}{2})=0$. Due to the map resulting in the same point regardless of $\theta_m$ for $j<k$, this results in a surjective mapping for all elements of the following form:
\end{flushleft}

\begin{center}
$\mathbb{B}_n\equiv \{\displaystyle\prod_{k=1}^{n}e^{\displaystyle i_k \theta_{k}}|$  $\theta_k=\displaystyle\frac{\pi}{2}, 1<k\leq n\}$
\end{center}

\begin{flushleft}
With a closed form expression for the elements that belong to the kernal of the map, we can explicitly define the inverse mapping for the isomorphism between  $\mathbb{S}^{2^{n}-1}_{\mathbb{G}_{n}}\setminus \mathbb{B}_n$ and $S^n$:
\end{flushleft}
\newcommand*\widefbox[1]{\fbox{\hspace{2em}#1\hspace{2em}}}

\fbox{
\parbox{\textwidth}{\centering
   $\mathbf{P}^{-1}\big(\bigg(\displaystyle\prod_{k=1}^{n} \cos \left(\theta_{k}\right)\bigg)e_0+\sum_{k=1}^{n}\bigg( \sin \left(\theta_{k}\right) \prod_{\ell=k+1}^{n} \cos \left(\theta_{\ell}\right) \prod_{m=1}^{k-1} (-1)^{\frac{(\theta_{m}-(\theta_{m} m o d \pi))}{\pi}}\bigg)e_{k}\big)
   \newline     \hphantom \quad\quad\quad\quad\quad\quad=\displaystyle\prod_{k=1}^{n}e^{\displaystyle i_k \theta_{k}}\newline  
  \theta_{m} \neq \frac{\pi}{2}, 1<m\leq n\quad\quad\quad$}}

\begin{flushleft}
Similarly, we can explicitly group all the elements in $\mathbb{B}_{n}$ that map to the same element in $S^n$:
\end{flushleft}
\begin{center}
$\mathbf{P}(\displaystyle\prod_{k=1}^{m-1}e^{\displaystyle i_k \theta_{k}}\cdot e^{\displaystyle \frac{i_m\pi}{2}}\cdot \displaystyle\prod_{\ell=m+1}^{n}e^{\displaystyle i_{\ell} \theta_{\ell}})=\mathbf{P}(\displaystyle e^{\displaystyle \frac{i_m\pi}{2}}\cdot \displaystyle\prod_{\ell=m+1}^{n}e^{\displaystyle i_{\ell} \theta_{\ell}})
$
\end{center}
\begin{center}
$\forall$  $\theta_m=\frac{\pi}{2},$ $1<m\leq n$
\end{center}
\begin{center}
\subsection{Commutative Diagram}
\end{center}
\begin{flushleft}
In order to generalize the surjective mapping from a partial torus onto the entire $S^n$, we must explicitly define our construction of a partial torus. In $\mathbb{R}^2$, the surfaces $\mathbb{T}^1$ and $S^1$ are equivalent, which means the construction begins to vary when $n>1$. Similar to the construction of polyspherical coordinates, we can generate $\mathbb{T}^n_{2^{1-n}}$ by scaling $\mathbb{T}^{n-1}_{2^{2-n}}$ with $\Big(a_n+cos(\theta_k)\Big)$ and adding an independent $(-1)^{\frac{(\theta_1-(\theta_1 m o d \pi))}{\pi}}\sin(\theta_k)$:
\end{flushleft}
\begin{definition}[Partial Torus in $\mathbb{R}^n$]
Let $a_m\geq1,\theta_1\in [0,2\pi),\theta_m\in [0,\pi),
\\1<j\leq n-1.$ We define the partial Torus as:
\end{definition}
\begin{center}
$\mathbb{T}^{n-1}_{2^{2-n}}\equiv\big\{\displaystyle\bigg(\displaystyle\cos \left(\theta_{1}\right)\prod_{k=2}^{n-1} (a_k+\cos \left(\theta_{k}\right))\bigg)e_0+\bigg(\displaystyle\sin \left(\theta_{1}\right)\prod_{k=2}^{n-1}(a_k+ \cos \left(\theta_{k}\right))\bigg)e_1
$
\end{center}
\begin{center}
$
+\displaystyle (-1)^{\frac{(\theta_1-(\theta_1 m o d \pi))}{\pi}} \sum_{k=2}^{n-1}\bigg(\sin \left(\theta_{k}\right) \prod_{\ell=k+1}^{n-1}(a_{\ell} + \cos \left(\theta_{\ell}\right)) \bigg)e_k\big\}$  
\end{center}
\begin{center}
We can observe how $\mathbb{T}^{n-1}_{2^{2-n}}$ generalizes to $\mathbb{T}^{n}_{2^{1-n}}$ using the following construction:
\end{center}
\begin{center}
$\mathbb{T}^{n}_{2^{1-n}}=\big\{(\mathbb{T}^{n-1}_{2^{2-n}}\big(a_k+ \cos(\theta_{k})\big),\sin(\theta_{k}))\big\}
$
\end{center}
\begin{center}
$
=\big\{\displaystyle\bigg(\Big(\displaystyle\cos \left(\theta_{1}\right)\prod_{k=2}^{n-1} (a_k+\cos \left(\theta_{k}\right))\Big)e_0+\Big(\displaystyle\sin \left(\theta_{1}\right)\prod_{k=2}^{n-1}(a_k+ \cos \left(\theta_{k}\right))\Big)e_1
$
\end{center}
\begin{center}
$
+\displaystyle (-1)^{\frac{(\theta_1-(\theta_1 m o d \pi))}{\pi}} \sum_{k=2}^{n-1}\Big(\sin \left(\theta_{k}\right) \prod_{\ell=k+1}^{n-1}(a_{\ell} + \cos \left(\theta_{\ell}\right)) \Big)e_k\bigg)\Big(a_n+\cos(\theta_n)\Big)
$
\end{center}
\begin{center}
$
+(-1)^{\frac{(\theta_1-(\theta_1 m o d \pi))}{\pi}}\sin(\theta_n)e_n|$  $a_m\geq1,\theta_1\in [0,2\pi),\theta_k\in [0,\pi), 1<k\leq n\big\} $
\end{center}
\begin{center}
$=\big\{\displaystyle\bigg(\displaystyle\cos \left(\theta_{1}\right)\prod_{k=2}^{n} (a_k+\cos \left(\theta_{k}\right))\bigg)e_0+\bigg(\displaystyle\sin \left(\theta_{1}\right)\prod_{k=2}^{n}(a_k+ \cos \left(\theta_{k}\right))\bigg)e_1
$
\end{center}
\begin{center}
$
+\displaystyle (-1)^{\frac{(\theta_1-(\theta_1 m o d \pi))}{\pi}} \sum_{k=2}^{n}\bigg(\sin \left(\theta_{k}\right) \prod_{\ell=k+1}^{n}(a_{\ell} + \cos \left(\theta_{\ell}\right)) \bigg)e_k|
$
\end{center}
\begin{center}
$
a_k\geq1,\theta_1\in [0,2\pi),\theta_k \in [0,\pi), 1<k\leq n\big\} $
\end{center}
\begin{flushleft}
Recall in Section 3.3, we gave the definition for the Map from $\mathbb{S}^{3}_{\mathbb{G}_{2}}$ to the Half Torus, which is a bijective mapping between the two spaces. The map can be defined as the following:
\begin{center}
\end{center}
\end{flushleft}
$\mathbf{Definition \; 3.5}$ (Map from $\mathbb{S}^3_{\mathbb{G}_2}$ to Half Torus). Let $\theta_1\in[0,2\pi), \theta_2\in[0,\pi) $
\begin{center}
$h(e^{\displaystyle i\theta_1}\cdot e^{\displaystyle j\theta_2})=\Big(\cos(\theta_1)\big(a_2+\cos(\theta_2)\big),\sin(\theta_1)\big(a_2+\cos(\theta_2)\big),(-1)^{\frac{(\theta_1-(\theta_1 m o d \pi))}{\pi}}\sin(\theta_2)\Big)$
\end{center}
\begin{center}
$h(\mathbb{S}^3_{\mathbb{G}_2})=\mathbb{T}^2_{2^{-1}}$
\end{center}
\begin{flushleft}
It can be trivially verified that all elements in $\mathbb{T}^{n}_{2^{1-n}}$ have a unique representation for all $\theta_1\in [0,2\pi)$ and $\theta\in [0,\pi), 1<j\leq n$. This allows us to generalize Definition 3.5 to a bijective mapping between $\mathbb{S}^{2^{n}-1}_{\mathbb{G}_{n}}$ and $\mathbb{T}^{n}_{2^{1-n}}$:
\end{flushleft}
\begin{definition}[Map from $\mathbb{S}^{2^n-1}_{\mathbb{G}_n}$ to $\mathbb{T}^{n}_{2^{1-n}}$] Let $a_k \geq 1,$ $\theta_1\in [0,2\pi),\\ \theta_k\in [0,\pi), 1< k\leq n$
\end{definition}
\begin{center}
$h(\displaystyle\prod^{n}_{k=1}e^{\displaystyle i_k\theta_k})= \bigg(\displaystyle\cos \left(\theta_{1}\right)\prod_{k=2}^{n} (a_k+\cos \left(\theta_{k}\right))\bigg)e_0+\bigg(\displaystyle\sin \left(\theta_{1}\right)\prod_{k=2}^{n}(a_k+ \cos \left(\theta_{k}\right))\bigg)e_1
$
\end{center}
\begin{center}
$
+\;\displaystyle (-1)^{\frac{(\theta_1-(\theta_1 m o d \pi))}{\pi}} \sum_{k=2}^{n}\bigg(\sin \left(\theta_{k}\right) \prod_{\ell=k+1}^{n}(a_{\ell} + \cos \left(\theta_{\ell}\right)) \bigg)e_k$
\end{center}
\begin{flushleft}
The final surjective morphism from $\mathbb{T}^{n}_{2^{1-n}}$ to $S^n$ can be defined as the limit of $\mathbb{T}^{n}_{2^{1-n}}$ as $a_k\rightarrow 0$ for all $k$:
\end{flushleft}
\begin{definition}[Map from $\mathbb{T}^{n}_{2^{1-n}}$ to $S^{n}$] Let $\theta_1\in [0,2\pi), \theta_k\in [0,\pi),\\ 1< k\leq n$
\end{definition}
\begin{center}
$\mu\bigg(\Big(\displaystyle\cos \left(\theta_{1}\right)\prod_{k=2}^{n} \big(a_k+\cos \left(\theta_{k}\right)\big)\Big)e_0+\Big(\displaystyle\sin \left(\theta_{1}\right)\prod_{k=2}^{n}\big(a_k+ \cos \left(\theta_{k}\right)\big)\Big)e_1
$
\end{center}
\begin{center}
$
+\;\displaystyle (-1)^{\frac{(\theta_1-(\theta_1 m o d \pi))}{\pi}} \sum_{k=2}^{n}\Big(\sin \left(\theta_{k}\right) \prod_{\ell=k+1}^{n}\big(a_{\ell} + \cos \left(\theta_{\ell}\right)\big) \Big)e_k\bigg)
$
\end{center}
\begin{center}
$
=\bigg(\displaystyle\prod_{k=1}^{n} \cos \left(\theta_{k}\right)\bigg)e_0+\bigg(\displaystyle\sin \left(\theta_{1}\right)\prod_{k=2}^{n} \cos \left(\theta_{k}\right)\bigg)e_1+\displaystyle (-1)^{\frac{(\theta_1-(\theta_1 m o d \pi))}{\pi}} \sum_{k=2}^{n}\bigg(\sin \left(\theta_{k}\right) \prod_{\ell=k+1}^{n} \cos \left(\theta_{\ell}\right) \bigg)e_k
$
\end{center}
\begin{center}
$
= \mathbf{P}(\displaystyle\prod^{n}_{k=1}e^{\displaystyle i_k\theta_k})$
\end{center}
\begin{center}
lim $a_k\rightarrow 0,$  $1<k\leq n $
\end{center}
\hfill
\begin{flushleft}
\newpage
The connection between $S_p^{2^{n}-1}$, $\mathbb{S}^{2^{n}-1}_{\mathbb{G}_{n}}$, $\mathbb{T}^{n}_{2^{1-n}}$, and $\mathbb{S}^n$ can be expressed through the following diagram:
\end{flushleft}

\[
\begin{tikzcd}
S_p^{2^{n}-1} 
\arrow[r,twoheadrightarrow,"f"]&
\mathbb{S}^{2^{n}-1}_{\mathbb{G}_{n}} \arrow[r,leftrightarrow,"h"]
\arrow[dr,swap,twoheadrightarrow,"\mathbf{P}"]
&
\mathbb{T}^n_{2^{1-n}}
\arrow[d,twoheadrightarrow,"\mu"]
\\
&&S^n
\end{tikzcd}
\]
\begin{center}
LG$(S_p^{2^{n}-1})=\mathbf{P}(f(S_p^{2^{n}-1}))=\mu(h(f(S_p^{2^{n}-1})))=S^n$
\end{center}
\begin{flushleft}
\section{Applications of the LG Fibration to Optimal Projections for Machine Learning}
\end{flushleft}
\begin{flushleft}
We can now look into properties such as the distance between two vectors before and after the transformation. In order to do this, we must define the inner product of two Multicomplex elements:
\end{flushleft}
\smallskip
\begin{center}
$\forall w_1,w_2\in\mathbb{C}_{n}: w_1=\displaystyle\prod^{n}_{k=1}e^{i_k\alpha_k},w_2=\displaystyle\prod^{n}_{k=1}e^{i_k\beta_k}$
\end{center}
\begin{center}
$\langle w_1,w_2\rangle=\langle \displaystyle\prod^{n}_{k=1}e^{i_k\alpha_k},\displaystyle\prod^{n}_{k=1}e^{i_k\beta_k}\rangle$
\end{center}
\begin{center}
$=\langle \displaystyle\prod^{n}_{k=1}e^{i_k(\alpha_k-\beta_k)},1\rangle$
\end{center}
\begin{center}
$=\langle \displaystyle\prod_{k=1}^{n} \cos (\alpha_k-\beta_k)+\sum_{k=1}^{n}i_k\big(\sin (\alpha_k-\beta_k) \prod_{\ell=k+1}^{n} \cos (\alpha_{\ell}-\beta_{\ell}) \prod_{m=1}^{k-1} e^{i_m (\alpha_{m}-\beta_{m})}\big),1\rangle$
\end{center}
\begin{center}
$=\langle \displaystyle\big(\prod_{k=1}^{n}\cos (\alpha_k-\beta_k)\big)e_0
$
\end{center}
\begin{center}
$
\displaystyle+\sum_{k=1}^{n}\Big(\sin (\alpha_k-\beta_k) \prod_{\ell=k+1}^{n} \cos (\alpha_{\ell}-\beta_{\ell}) \prod_{m=1}^{k-1} e^{i_m (\alpha_{m}-\beta_{m})}\Big)e_k,e_0\rangle$
\end{center}
\begin{center}
$=\displaystyle\prod^{n}_{k=1}\cos(\alpha_k-\beta_k)$
\end{center}

\begin{flushleft}
This elegant algebraic result comes from observing that the real component of the product of rotations must be equal to the product of the real component from each independent rotation group. Furthermore, we use the equivalence that the inner product between elements in $\mathbb{C}_{n}$ can be interpreted as the standard inner product between their vector representations in $\mathbb{R}^{2^{n}}$\cite{luna2015bicomplex}. We now have the tools required to compare the inner product between two vectors before and after the LG Fibration Projection:
\end{flushleft}
\begin{center}
$ w_1=e^{\displaystyle i_1\alpha_1}\cdot e^{\displaystyle i_2\alpha_2} $
\end{center}
\begin{center}
$\mathbf{P}(w_1)=\cos(\alpha_1)\cos(\alpha_2)e_0+\sin(\alpha_1)\cos(\alpha_2)e_1+(-1)^{\frac{(\alpha_1-(\alpha_1 m o d \pi))}{\pi}}\sin(\alpha_2)e_2$
\end{center}
\begin{center}
$w_2= e^{\displaystyle i_1\beta_1}\cdot e^{\displaystyle i_2\beta_2} $
\end{center}
\begin{center}
$\mathbf{P}(w_2)=\cos(\beta_1)\cos(\beta_2)e_0+\sin(\beta_1)\cos(\beta_2)e_2+(-1)^{\frac{(\beta_1-(\beta_1 m o d \pi))}{\pi}}\sin(\beta_2)e_2$
\end{center}
\
\begin{flushleft}
We can define a Distance Difference function to measure the difference between the inner product of the vectors before and after the modular transformation:
\end{flushleft}
\smallskip
\begin{center}
$D(\alpha_1,\alpha_2,\beta_1,\beta_2)\equiv\big|\langle w_1,w_2\rangle-\langle \mathbf{P}(w_1),\mathbf{P}(w_2)\rangle\big|$
\end{center}
\begin{center}
$=\Big|\cos(\alpha_1-\beta_1)\cos(\alpha_2-\beta_2)-\Big(\cos(\alpha_1)\cos(\alpha_2)\cos(\beta_1)\cos(\beta_2)+\sin(\alpha_1)\cos(\alpha_2)\sin(\beta_1)\cos(\beta_2)+(-1)^{\frac{(\alpha_1-\alpha_1 m o d \pi)}{\pi}}(-1)^{\frac{(\beta_1-\beta_1 m o d \pi)}{\pi}}\sin(\alpha_2)\sin(\beta_2)\Big)\Big|$
\end{center}
\begin{center}
$=\Big|\cos(\alpha_1-\beta_1)\cos(\alpha_2-\beta_2)-\Big(\cos(\alpha_2)\cos(\beta_2)\big(\cos(\alpha_1)\cos(\beta_1)+\sin(\alpha_1)\sin(\beta_1)\big)+(-1)^{\frac{(\alpha_1-\alpha_1 m o d \pi)}{\pi}}(-1)^{\frac{(\beta_1-\beta_1 m o d \pi)}{\pi}}\sin(\alpha_2)\sin(\beta_2)\Big)\Big|$
\end{center}
\begin{center}
$=\Big|\cos(\alpha_1-\beta_1)\cos(\alpha_2-\beta_2)-\cos(\alpha_2)\cos(\beta_2)\cos(\alpha_1-\beta_1)-(-1)^{\frac{(\alpha_1-\alpha_1 m o d \pi)}{\pi}}(-1)^{\frac{(\beta_1-\beta_1 m o d \pi)}{\pi}}\sin(\alpha_2)\sin(\beta_2)\Big|$
\end{center}
\begin{center}
$=\Big|\cos(\alpha_1-\beta_1)\Big(\cos(\alpha_2-\beta_2)-\cos(\alpha_2)\cos(\beta_2)\Big)-(-1)^{\frac{(\alpha_1-\alpha_1 m o d \pi)}{\pi}}(-1)^{\frac{(\beta_1-\beta_1 m o d \pi)}{\pi}}\sin(\alpha_2)\sin(\beta_2)\Big|$
\end{center}
\begin{center}
$=\Big|\cos(\alpha_1-\beta_1)\sin(\alpha_2)\sin(\beta_2)-(-1)^{\frac{(\alpha_1+\beta_1-\alpha_1 m o d \pi-\beta_1 m o d \pi)}{\pi}}\sin(\alpha_2)\sin(\beta_2)\Big|$
\end{center}
\begin{center}
$=\Big|\sin(\alpha_2)\sin(\beta_2)\big(\cos(\alpha_1-\beta_1)-(-1)^{\frac{(\alpha_1+\beta_1-\alpha_1 m o d \pi-\beta_1 m o d \pi)}{\pi}}\big)\Big|$
\end{center}
\smallskip
\begin{flushleft}
We can algebraically determine that the inner product between two vectors before and after the transformation is invariant if:
\end{flushleft}

\begin{center}
$\displaystyle \sin(\alpha_2)\sin(\beta_2)=0$ or $\cos(\alpha_1-\beta_1)=(-1)^{\frac{(\alpha_1+\beta_1-\alpha_1 m o d \pi-\beta_1 m o d \pi)}{\pi}}$
\end{center}
\begin{flushleft}
The General Difference function can be defined as the following:
\end{flushleft}

\begin{center}
$ \boxed{D(\alpha_{1},...,\alpha_{n},\beta_{1},...,\beta_{n})=\big|\displaystyle\prod^n_{k=1}\cos(\alpha_k-\beta_k)-\langle \mathbf{P}(\displaystyle\prod_{k=1}^{n}e^{\displaystyle i_k \alpha_{k}}),\mathbf{P}(\displaystyle\prod_{k=1}^{n}e^{\displaystyle i_k \beta_{k}})\rangle\big|}$
\end{center}

\begin{flushleft}
Minimizing the function determines which pairs of vectors have an inner product that is invariant under the LG Fibration. This means that the distance between a pair of vectors in the original space as well as in the dimensionally reduced space is equivalent, which makes the General Difference function equal to zero. Computationally solving for these regions allows for the optimal distribution of data for minimal information loss from applying the LG Fibration.
\end{flushleft}
\begin{flushleft}
\section{Selected Graphs}
\end{flushleft}

We can observe parameterizations of $\mathbb{S}^{3}_{\mathbb{G}_2}$ projected down onto $S^2$, and compare how the projection differs from a standard parameterization of polyspherical Coordinates for $S^2$. The LG Fibration Projection $\mathbf{P}(e^{\displaystyle i_1\theta_1}\cdot e^{\displaystyle i_2\theta_2})$ takes the following form:
\begin{center}
$x'=\cos(\theta_1)\cos(\theta_2)$
\end{center}
\begin{center}
$y'=\sin(\theta_1)\cos(\theta_2)$
\end{center}
\begin{center}
$z'=(-1)^{\frac{(\theta_{1}-(\theta_{1} m o d \pi))}{\pi}}\sin(\theta_2)$
\end{center}
\begin{flushleft}
We can create a parametric line of a projected rotation through the following constraint:
\end{flushleft}
\begin{center}
$\theta_2=a\theta_1$
\end{center}
\begin{center}
$x'=\cos(\theta_1)\cos(a\theta_1)$
\end{center}
\begin{center}
$y'=\sin(\theta_1)\cos(a\theta_1)$
\end{center}
\begin{center}
$z'=(-1)^{\frac{(\theta_{1}-(\theta_{1} m o d \pi))}{\pi}}\sin(a\theta_1)$
\end{center}
\begin{center}
The following examples allow us to gain heuristic understanding of how the periodic behavior of the higher dimensional rotations is preserved by the map:
\end{center}
\newpage
\begin{center}
$\theta_2=\theta_1$
\end{center}
\begin{tabular}{|c|c|}
      \hline
      \includegraphics[width=0.5\textwidth]{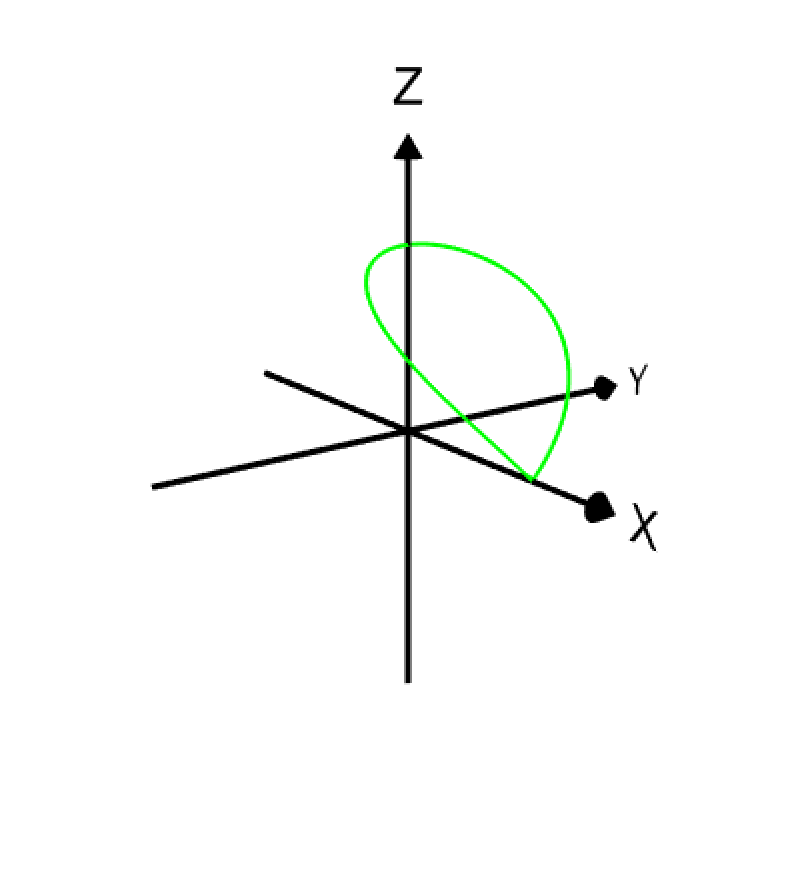}&
      \includegraphics[width=0.5\textwidth]{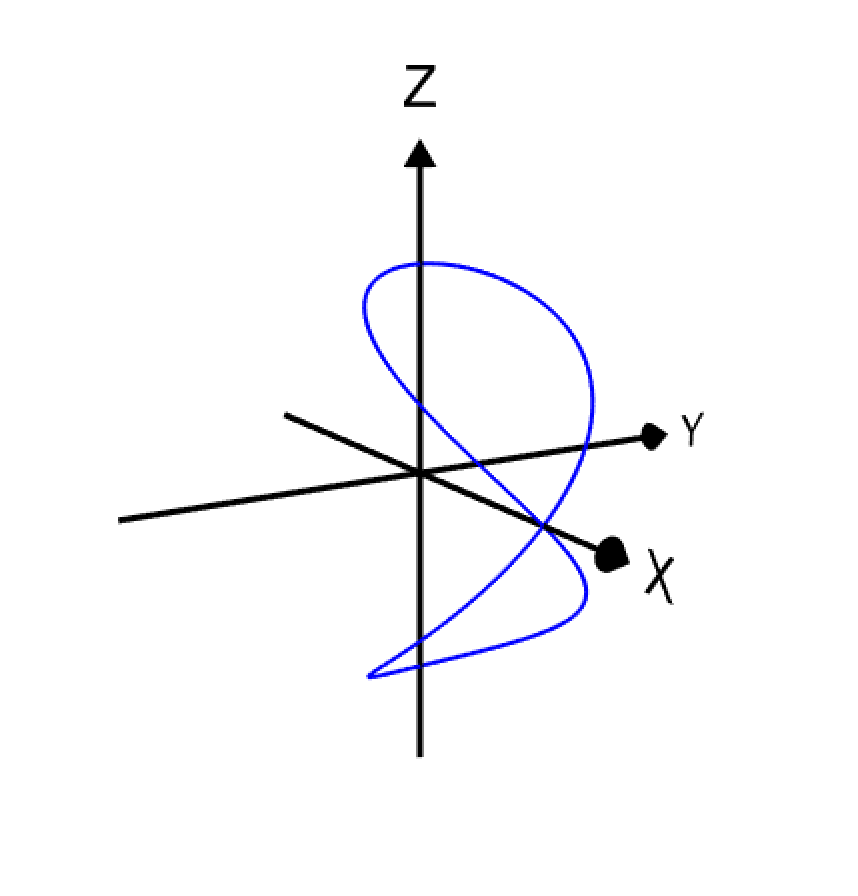}
\\
      $x'=\cos(\theta_1)\cos(\theta_1)$&$x=\cos(\theta_1)\cos(\theta_1)$ \\
       $y'=\sin(\theta_1)\cos(\theta_1)$ & $y=\sin(\theta_1)\cos(\theta_1)$ 
       \\
       $z'=(-1)^{\frac{(\theta_{1}-(\theta_{1} m o d \pi))}{\pi}}\sin(\theta_1)$ & $z=\sin(\theta_1)$
       \\
   \hline
\end{tabular}
\begin{center}
$\theta_2=2\theta_1$
\end{center}
\begin{center}
\begin{tabular}{|c|c|}
      \hline
      \includegraphics[width=0.5\textwidth]{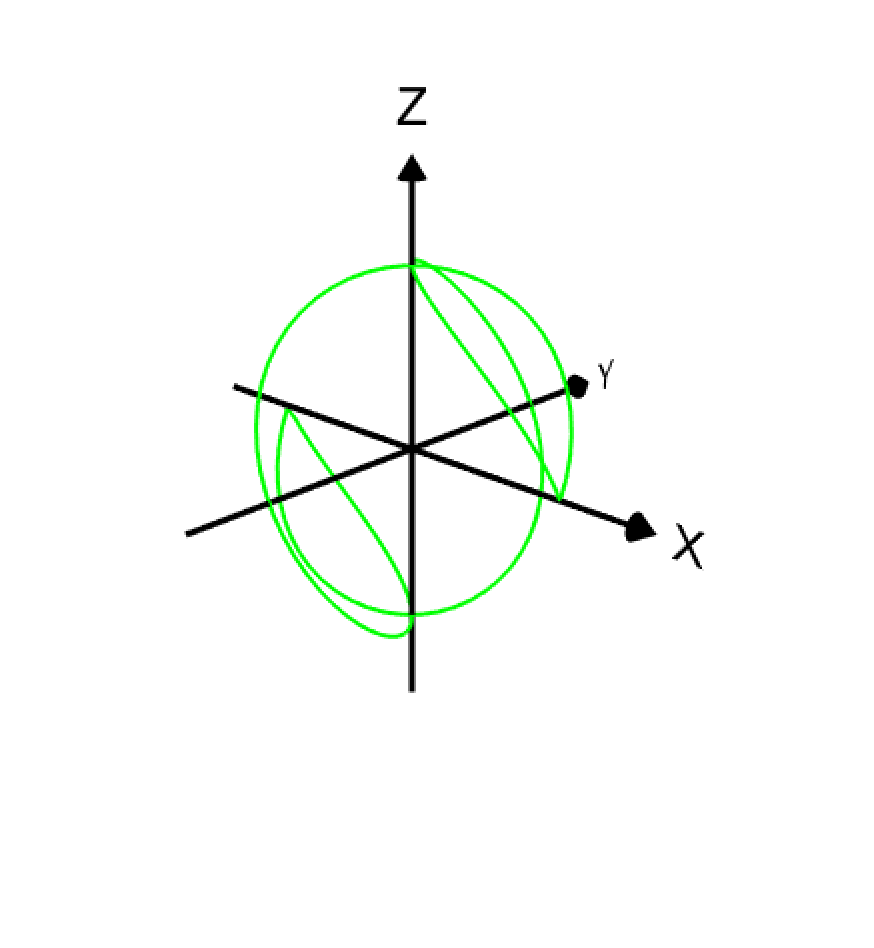}&
     \includegraphics[width=0.5\textwidth]{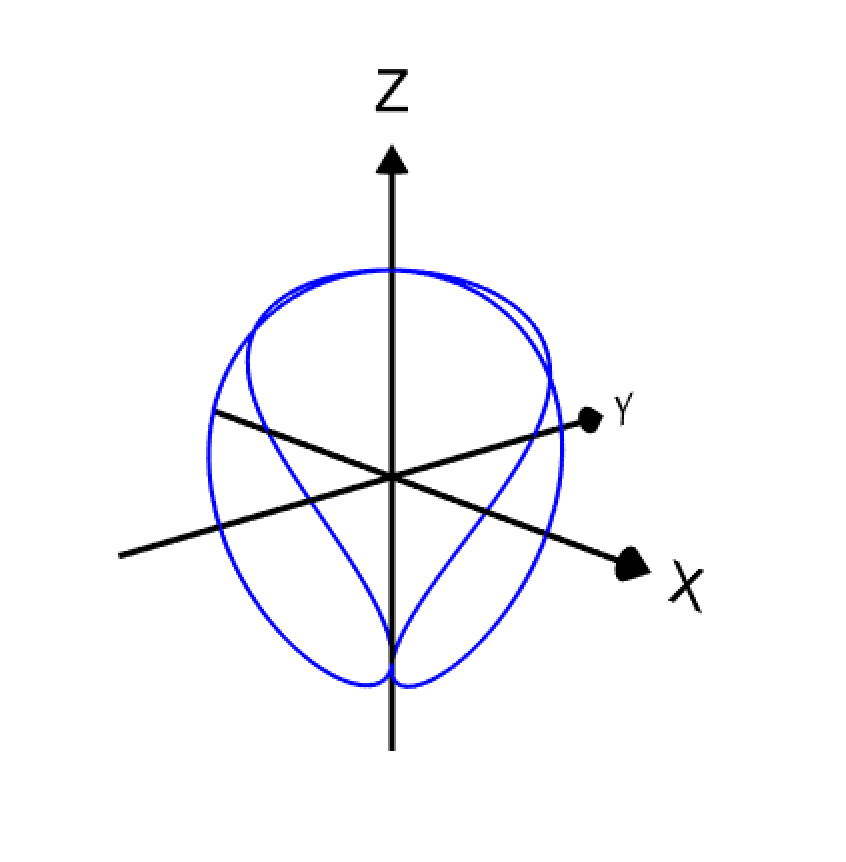}
\\
      $x'=\cos(\theta_1)\cos(2\theta_1)$&$x=\cos(\theta_1)\cos(2\theta_1)$ \\
       $y'=\sin(\theta_1)\cos(2\theta_1)$ & $y=\sin(\theta_1)\cos(2\theta_1)$ 
       \\
       $z'=(-1)^{\frac{(\theta_{1}-(\theta_{1} m o d \pi))}{\pi}}\sin(2\theta_1)$ & $z=\sin(2\theta_1)$
       \\
   \hline
\end{tabular}
\end{center}
\smallskip
\begin{center}
Overlapping the LG Fibration with the polyspherical graph allows us to see the effect of the modular shift. We can also observe that the mapping produces 1 non-differentiable point if $a$ is odd, and 2 non-differentiable points if $a$ is even.
\end{center}
\begin{center}
\begin{tabular}{|c|c|}
      \hline
      \includegraphics[width=0.5\textwidth]{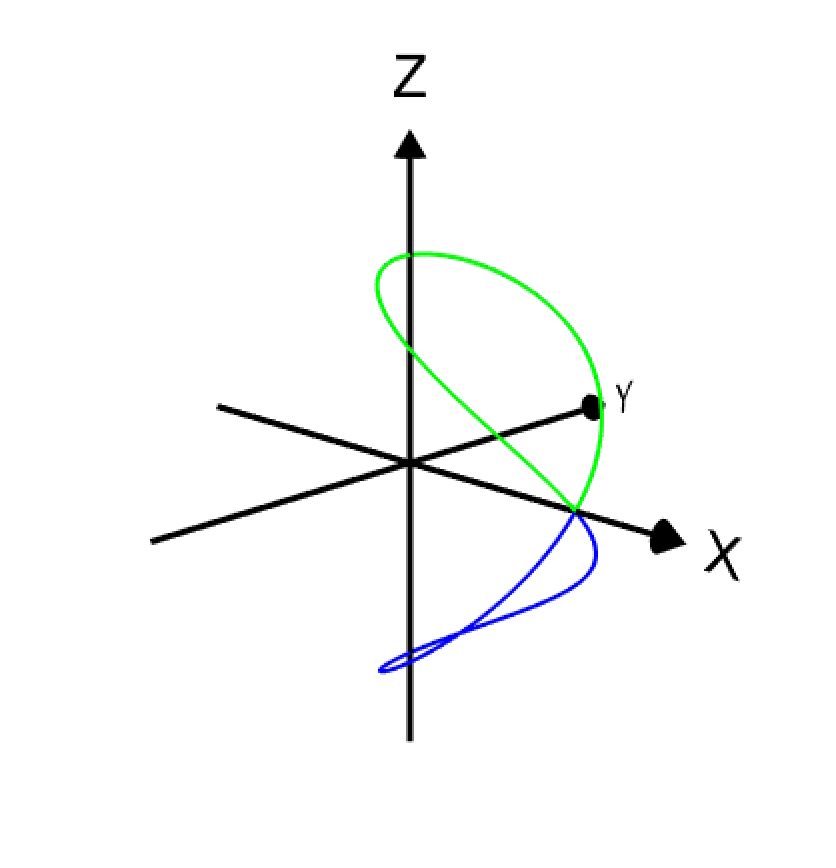}&
      \includegraphics[width=0.5\textwidth]{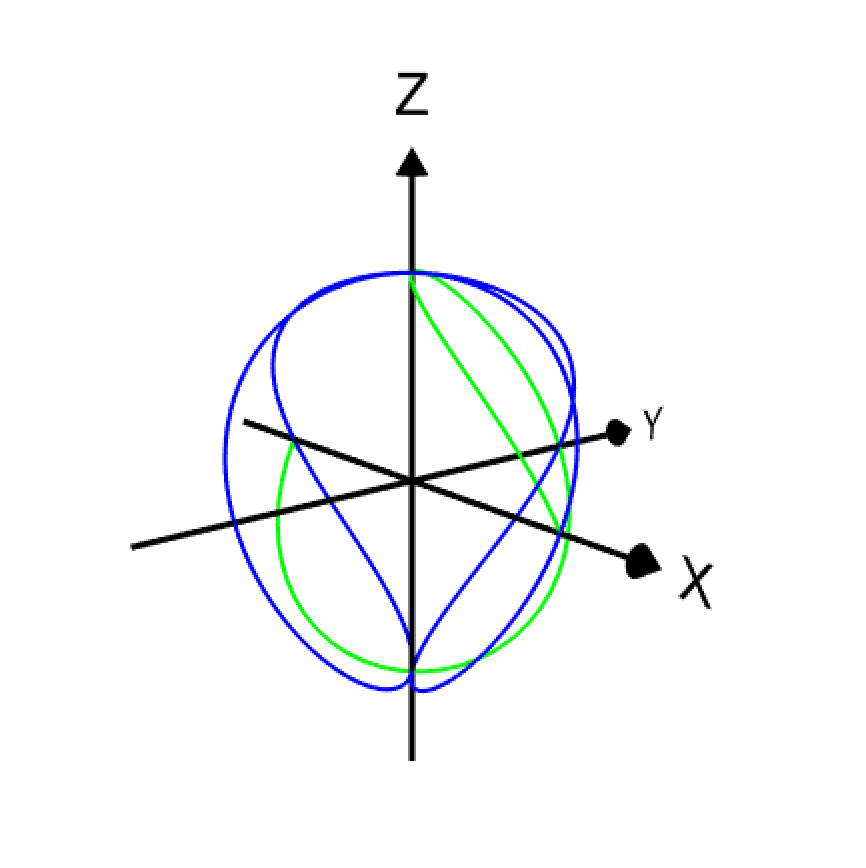}
      \\
   \hline
\end{tabular}
\end{center}
\begin{center}
\end{center}

\begin{center}
Finally, the XY cross sections form rose curves with $2a$ rose petals if $a$ is even, and $a$ rose petals if $a$ is odd:
\end{center}
\begin{center}
\begin{tabular}{|c|c|}
      \hline
      \includegraphics[width=0.5\textwidth]{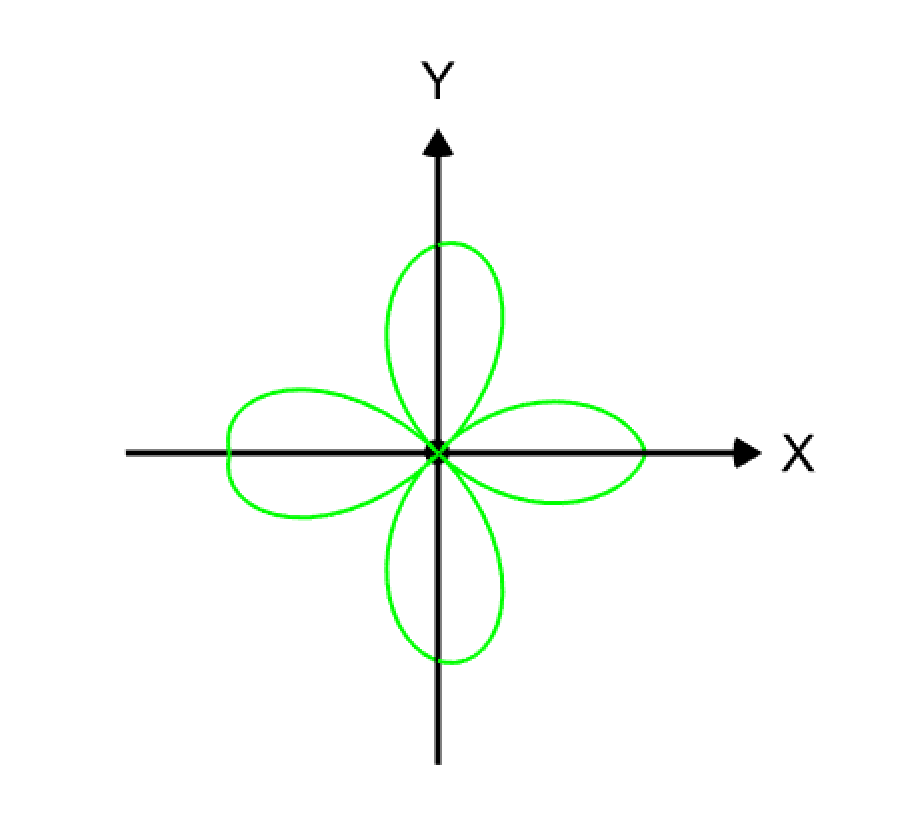}&
      \includegraphics[width=0.5\textwidth]{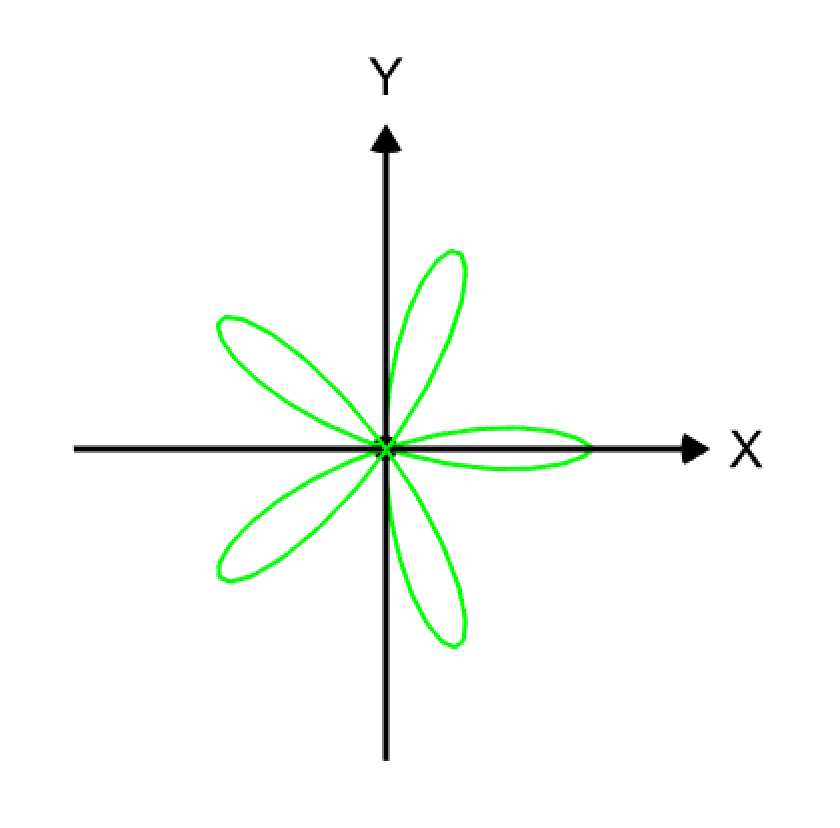}
\\
      $\theta_2=2\theta_1$&$\theta_2=5\theta_1$ \\
    
   \hline
\end{tabular}
\end{center}
\begin{flushleft}
\section{Conclusion}
\end{flushleft}
\begin{flushleft}
\quad Our investigation began with the derivation of polyspherical coordinates, which we used to bridge a link with a subgroup of the general Abelian Multicomplex Rotation Group. The bijective relationship between $\mathbb{S}^{2^{n}-1}_{\mathbb{G}_{n}}$ and $\mathbb{T}^{n}_{2^{1-n}}$ allows for the LG fibration to map surjectively onto $S^n$ from $S^{2^n-1}$. This procedure immediately opens up the possibility for two separate areas of research: $\mathbb{R}^{n+1}$ rotations using commutative $\mathbb{R}^{2^{n}}$ dimensional subgroups in pure mathematics, and Dimensionality reduction algorithms for sparse high dimensional spaces in theoretical computer science.
\end{flushleft}
\begin{flushleft}
\quad Large scale implementation of the LG Fibration in industrial Dimensionality Reduction requires further research into methods of preconditioning data sets such that the inital data will optimally lie within the defined manifold. Additional research may also be conducted into determining other suitable low dimensional manifolds that may be embedded in high dimensional space, for this would allow alternative approaches to minimize data loss. Numerous graphing and visualization topics may build off of this research as well, for the long history of Quaternions being used to model 3D rotations may be rivaled by the appropriate application of Bicomplex Numbers.\cite{vicci2001quaternions}
\end{flushleft}
\begin{flushleft}
\quad This paper utilizes a variety of algebraic and geometric properties present in the Multicomplex numbers to derive its results, however numerous questions regarding the general group theoretic symmetries of the structure are yet to be resolved and thus require extensive investigation. Future publications should aim to utilize the abelian rotation and hyperbolic groups in $\mathbb{C}_n$ for mathematical physics, such as generalizing how the split complex numbers may be used to perform a Lorentz boost in 1+1 Minkowski space.\cite{catoni2008mathematics} Discovering the LG Fibration was made possible because of the clear connection between geometry and algebra in the Multicomplex numbers, and our hope through this paper is to inspire further research into the applications of geometric groups embedded in high dimensional spaces.
\end{flushleft}
\begin{flushleft}
 \printbibliography
\end{flushleft}
\end{document}